\documentclass[11pt]{article}

\usepackage[utf8]{inputenc}
\usepackage{amsmath}
\usepackage{amsfonts,amssymb,latexsym,multicol}
\usepackage[francais]{babel}
\usepackage[T1]{fontenc}

\newcounter{fig}
\newtheorem{theo}{Th\'eor\`eme}
\newtheorem{lem}{Lemme}

\newtheorem{prop}{Proposition}

\newcommand{\cad}{\text{c'est-\`a-dire }}
\newcommand{\expli}[1]{\quad\text{\footnotesize (#1)}}

\makeatletter
\ifcase\@ptsize
\def\one{{\mathchoice{\rm 1\mskip-4.3mu l}{\rm 1\mskip-4.3mu l}
                     {\rm 1\mskip-4.5mu l}{\rm 1\mskip-5.0mu l}}}
\or
\def\one{{\mathchoice{\rm 1\mskip-4.7mu l}{\rm 1\mskip-4.7mu l}
                     {\rm 1\mskip-4.7mu l}{\rm 1\mskip-5.0mu l}}}
\or
\def\one{{\mathchoice{\rm 1\mskip-4.4mu l}{\rm 1\mskip-4.4mu l}
                     {\rm 1\mskip-4.6mu l}{\rm 1\mskip-5.0mu l}}}
\fi
\makeatother

\newcommand{\eps}{\varepsilon}
\newcommand{\fhi}{\varphi}
\newcommand{\ioe}{\leqslant}
\newcommand{\soe}{\geqslant}
\renewcommand{\le}{\leqslant}
\renewcommand{\ge}{\geqslant}
\newcommand{\vers}{\rightarrow}

\newcommand{\demi}{{\frac{1}{2}}}

\newcommand{\ci}{{\rm ci}}
\newcommand{\Ci}{{\rm Ci}}

\newcommand{\Wcal}{{\mathcal W}}

\newcommand{\Cgot}{{\mathfrak C}}

\newcommand{\Fgot}{{\mathfrak F}}

\newcommand{\cgot}{{\mathfrak c}}

\newcommand{\Nat}{{\mathbb N}}
\newcommand{\Int}{{\mathbb Z}}
\newcommand{\Rat}{{\mathbb Q}}
\newcommand{\Real}{{\mathbb R}}

\newcommand{\sgn}{{\rm sgn}}

\newcommand{\fin}{\hfill$\Box$}
\newcommand{\dem}{\noindent {\bf D\'emonstration\ }}
\newcommand{\fine}{\tag*{\mbox{$\Box$}}}

\providecommand{\bysame}{\leavevmode ---\ }
\providecommand{\og}{``}
\providecommand{\fg}{''}
\providecommand{\smfandname}{et}
\providecommand{\smfedsname}{\'eds.}
\providecommand{\smfedname}{\'ed.}

\newcommand{\W}{\mathcal{W}}
\newcommand{\intW}{\Upsilon}

\title{Sur l'autocorr\'elation multiplicative\\ de la fonction {\og partie fractionnaire\fg}\\ et une fonction définie par J. R. Wilton}

\author{Michel Balazard et Bruno Martin}

\begin{document}
\maketitle

\begin{flushright}
\textit{\`A la mémoire de notre ami Patrick Sargos.}
\end{flushright}

\begin{center}
  {\sc Abstract}
\end{center}
\begin{quote}
{\footnotesize We describe the points of differentiability of the function
\begin{equation*}
 A(\lambda)= \int_0^\infty \{t \} \{\lambda t \}\frac{dt}{t^2}, 
\end{equation*}
where $\{t\}$ denotes the fractional part of $t$. In connection with this question, we study series involving the first Bernoulli function, the arithmetical function "number of divisors", and the Gauss map $\alpha(x)=\{1/x\}$ from the theory of continued fractions. A key role is played by a function defined in 1933 by J. R. Wilton, similar to the Brjuno function of dynamical systems theory. A unifying theme of our exposition is the use of functional equations involving the Gauss map, allowing us to reprove and refine a theorem of Wilton, la Bretèche and Tenenbaum.
}
\end{quote}

\begin{center}
  {\sc Keywords}
\end{center}
\begin{quote}
{\footnotesize Fractional part, Autocorrelation, Continued fractions, Approximate functional equations }
\end{quote}

\newpage

\tableofcontents

\section{Introduction}

Posons pour $\lambda \soe 0$
\begin{equation*}
 A(\lambda)= \int_0^\infty \{t \} \{\lambda t \}\frac{dt}{t^2}, 
\end{equation*}
o\`u $\{t\}=t-\lfloor t\rfloor$ d\'esigne la partie fractionnaire du nombre r\'eel
$t$, et $\lfloor t\rfloor$ sa partie enti\`ere.
La fonction $A(\lambda)$ est la \emph{fonction d'autocorr\'elation multiplicative de la
fonction {\og partie fractionnaire\fg}}, introduite  par B\'aez-Duarte {\it et al.}
(cf. \cite{baez-duarte-all}) dans le contexte de l'\'etude du crit\`ere de Nyman pour
l'hypoth\`ese de Riemann. 

B\'aez-Duarte {\it et al.} \'etablissent en particulier que la fonction $A$, continue sur $[0,\infty[$,  n'est
d\'erivable en aucun point rationnel et posent la question de d\'eterminer
l'ensemble de ses points de d\'erivabilit\'e. L'objet de la première partie du pr\'esent travail (\S\S\ref{par:fractions-continues}-\ref{t160}) est de
r\'epondre compl\`etement \`a cette question.

\begin{theo}\label{t114}
 Les points de d\'erivabilit\'e de la fonction $A$ sont les nombres irrationnels positifs $\lambda$ tels que la s\'erie
$$
\sum_{k\soe 0}(-1)^k\frac{\log q_{k+1}}{q_k}
$$
converge, o\`u les $q_k$ sont les d\'enominateurs des r\'eduites de $\lambda$. 
\end{theo}

Pour démontrer le théorème \ref{t114}, nous ramenons l'étude de la fonction d'autocorrélation $A$ à celle d'une fonction $\Wcal$, somme d'une série introduite par Wilton en 1933 (cf. \cite{Wilton}) dans l'étude de la série trigonométrique
\begin{equation}\label{t152}
\psi_1(x)=-\frac{1}{\pi}\sum_{n\soe 1}\frac{\tau(n)}{n}\sin(2\pi nx)\, ,
\end{equation}
o\`u $\tau(n)$ d\'esigne le nombre de diviseurs du nombre entier naturel $n$. La série $\W(x)$ est définie au \S\ref{t94} en termes du développement de $x$ en fraction continue, et nous verrons alors qu'elle définit une fonction intégrable sur $[0,1]$, solution de l'\'equation fonctionnelle
\begin{equation}
  \label{t93}
  \Wcal (x)=\log (1/x) -x\Wcal(\{1/x\}).
\end{equation}

Il s'av\`ere que la fonction d'autocorr\'elation $A$ et la primitive
$$
\Upsilon (x)=\int_0^x\Wcal(t)dt 
$$
de la fonction de Wilton $\Wcal$ sont intimement li\'ees.
\begin{prop}\label{t95}
Pour $0<\lambda <1$, on a
$$
A(\lambda) = \rho(\lambda)-\frac{\Upsilon(\lambda)}{2\lambda},
$$
o\`u $\rho$ est une fonction continue sur $]0,1[$, et d\'erivable en chaque irrationnel.
\end{prop}

Ainsi, la recherche des points de d\'erivabilit\'e de $A$ est ramen\'ee \`a celle des points de d\'erivabilit\'e de $\Upsilon$. En effet la fonction $A$ n'est pas d\'erivable \`a droite en $0$ : on a $A(\lambda)\sim\demi\lambda\log (1/\lambda)$ quand $\lambda \vers 0$ (relation \eqref{t202} \textit{infra}). D'autre part, l'identit\'e $A(\lambda)=\lambda A(1/\lambda)$ permet de restreindre l'\'etude \`a l'intervalle $]0,1[$. 

Cette question est tr\`es voisine de celle trait\'ee dans notre pr\'ec\'edent article \cite{Brjuno} concernant le comportement local moyen de la fonction de Brjuno, solution de l'\'equation fonctionnelle
\begin{equation}
  \label{t96}
  \Phi (x)=\log (1/x) +x\Phi(\{1/x\}),
\end{equation}
qui ne diff\`ere de \eqref{t93} que par un signe.

La fonction de Brjuno joue un r\^ole important dans la th\'eorie des syst\`emes dynamiques, plus pr\'ecis\'ement dans l'\'etude des it\'erations d'un polyn\^ome quadratique (cf. par exemple l'article \cite{BC} qui d\'emontre un r\'esultat substantiel dans ce domaine, et renvoie aux r\'ef\'erences ant\'erieures). Observons cependant qu'elle appara\^{\i}t d\'ej\`a en filigrane dans l'article de Wilton (cf. \cite{Wilton} (7.32), p. 235). 

Pour \'etudier $\Phi$ et notamment d\'eterminer l'ensemble de ses points de
 Lebesgue, c'est-\`a-dire, les points $x\in\Real$  tel que
\begin{equation}
 \lim_{h\to 0} \frac{1}{h} \int_{x}^{x+h}|\Phi(t)-\Phi(x)|dt=0, 
\end{equation}
nous avons d\'evelopp\'e dans \cite{Brjuno} plusieurs \'el\'ements th\'eoriques concernant
le comportement en moyenne de fonctions directement li\'ees au d\'eveloppement en
fraction continue d'un nombre r\'eel.
Dans le présent article, nous rappelons et exploitons ces r\'esultats pour \'etudier le
comportement local moyen de la fonction $\W$ de Wilton. Cela nous permet de montrer que les points de d\'erivabilit\'e de $\intW$, et donc ceux de $A$, sont pr\'ecis\'ement ceux mentionn\'es dans le th\'eor\`eme \ref{t114}. 

\smallskip

Dans la d\'emonstration de la proposition \ref{t95}, la fonction auxiliaire suivante joue un r\^ole essentiel :
\begin{equation}\label{eq:def-phiun}
 \fhi_1(t)= \sum_{n\ge 1} \frac{ B_1(nt)}{n},
\end{equation} 
o\`u $B_1$ d\'esigne la premi\`ere fonction de Bernoulli normalis\'ee :
\begin{equation*}
 B_1(t)= \{t\}-1/2 +[t\in \Int]/2,
\end{equation*}
où $[P]$ désigne la valeur de vérité de l'assertion $[P]$ (notation d'Iverson).

La s\'erie $\fhi_1(t)$ converge presque partout et d\'efinit une fonction p\'eriodique, de p\'eriode $1$, int\'egrable sur $[0,1]$. B\'aez-Duarte \textit{et al.} ont \'etabli dans \cite{baez-duarte-all} la relation suivante entre $A$ et $\fhi_1$ : pour tout $\lambda>0$, 
\begin{equation}\label{t151}
 A(\lambda)=\demi\log \lambda
+\frac{A(1)+1}{2}-\lambda\int_{\lambda}^{\infty}\fhi_1(t)\frac{dt}{t^2}. 
\end{equation}
Pour démontrer la proposition \ref{t95}, il reste donc à établir une relation entre les fonctions $\fhi_1$ et $\Wcal$. C'est l'objet de la proposition suivante, démontrée sous une forme plus précise au \S\ref{t161} (proposition \ref{prop:lien-fhi1-W}).
\begin{prop}\label{t107}
On a  
\begin{equation}\label{egalitepp}
\fhi_1(x)=-\frac12 \W(x)+G(x)\quad (p. p. ),
\end{equation}
o\`u $G:\Real\to\Real$ est une fonction $1$-p\'eriodique born\'ee, continue en tout irrationnel.  
 \end{prop}

La fonction $G$ est définie et étudiée au \S\ref{t163}. La proposition \ref{t95} r\'esulte alors de \eqref{t151} et \eqref{egalitepp} par une int\'egration par parties.

\smallskip

Nous avons choisi de rédiger l'ensemble des raisonnements menant au théorème \ref{t114} dans un cadre conceptuel ajusté à l'objet étudié et au résultat présenté. Alors que B\'aez-Duarte {\it et al.} montraient le parti que l'on pouvait tirer, pour l'\'etude de la fonction $A$, de la repr\'esentation de $\fhi_1$ par la s\'erie trigonom\'etrique $\psi_1$ définie par \eqref{t152}, notre  démonstration se passe enti\`erement de cette repr\'esentation\footnote{En cohérence avec ce parti pris, nous écrirons $A(1)$ et $\zeta(2)$ sans expliciter leurs valeurs respectives $\log(2\pi) -\gamma$ et $\pi^2/6$.}. Nous éliminons également le recours à la transformation de Mellin, utilisée dans \cite{baez-duarte-all} pour démontrer \eqref{t151}, et redémontrons cette identité par des manipulations élémentaires de séries et d'intégrales au \S\ref{t162}.

\medskip

Dans la deuxième partie du présent article (\S\S\ref{t120}-\ref{t149}), issue des recherches effectuées pour démontrer la proposition \ref{t107}, nous étudions et précisons un théorème de Wilton, la Bretèche et Tenenbaum. Nous aboutissons à l'énoncé suivant.

\begin{theo}\label{t150}
Les trois s\'eries $\fhi_1(x)$, $\psi_1(x)$ et $\W(x)$ convergent pour les m\^emes valeurs de la variable $x$. En tout point de convergence $x$, on a  
$$
\fhi_1(x)=\psi_1(x)=-\frac12 \W(x)+G(x)+\delta(x),
$$
o\`u $\delta:\Real\to\Real$ est  une fonction $1$-p\'eriodique nulle sur $\Real\setminus\Rat$.  
 \end{theo}
\medskip

L'assertion suivant laquelle les points de convergence de $\psi_1$ et de $\W$ sont les m\^emes est due \`a Wilton (1933, \cite{Wilton}, $(20_{III})$, p. 221). Ce r\'esultat fut retrouv\'e en 2004 par la Bret\`eche et Tenenbaum, qui d\'emontr\`erent \'egalement que les points de convergence de $\fhi_1$ co\"{\i}ncidaient avec ceux de $\psi_1$ et $\W$, ainsi que l'identit\'e $\fhi_1=\psi_1$ en tout point de convergence (\cite{breteche-tenenbaum}, Th\'eor\`eme 4.4, p. 16). Enfin, l'\'egalit\'e $\fhi_1=-\tfrac12\W+G+\delta$ est notre contribution \`a cet \'enonc\'e. Il est complété par la définition et les propriétés de la fonction $G$, qui sont l'objet du \S\ref{t163} ; quant à la fonction $\delta$, elle est définie par \eqref{t122} au \S\ref{t149}.

Cela \'etant, notre but dans cette seconde partie est surtout de fournir l'expos\'e autonome d'une d\'emonstration du th\'eor\`eme \ref{t150}, apportant à la question étudiée un éclairage complémentaire aux approches de Wilton et de La Bret\`eche et Tenenbaum, en la pla{\c c}ant, dans une perspective historique, au carrefour de plusieurs probl\'ematiques. En particulier, le leitmotiv technique de notre expos\'e est l'utilisation d'\'equations fonctionnelles du type
\begin{equation}\label{t115}
f(x)=-xf(\{1/x\})+g(x) \quad (x \in X),
\end {equation}
o\`u $g$ est une fonction suppos\'ee connue et $f$ la fonction inconnue, et d'\'equations fonctionnelles approch\'ees apparent\'ees \`a \eqref{t115}. C'\'etait d\'ej\`a le cas dans \cite{Wilton}, et Wilton cite comme source de cette id\'ee l'article \cite{HL} de Hardy et Littlewood. 

Ainsi la fonction de Wilton v\'erifie l'\'equation \eqref{t115} avec $g(x)=\log(1/x)$ (cf. \S\ref{t116}), et les sommes partielles de $\fhi_1$ et $\psi_1$ v\'erifient des \'equations fonctionnelles approch\'ees similaires \`a \eqref{t115} avec des fonctions $g(x)=-\demi\log(1/x)+O(1)$ (cf. \S\ref{par:fhi1} et \S\ref{par:psi1}). C'est l\`a le fait essentiel qui permet de d\'emontrer le th\'eor\`eme. La fonction $G$ appara\^{\i}t comme solution d'une \'equation fonctionnelle \eqref{t115} avec une fonction $g$ continue sur $[0,1]$, construite \`a partir de la fonction $A$. Celle-ci joue donc ici un r\^ole auxiliaire, alors qu'elle est la premi\`ere motivation de ce travail. 

\medskip

Voici quel est le plan de cet article. 

\emph{Première partie : points de dérivabilité de la fonction $A$.} Au \S\ref{par:fractions-continues} nous rappelons, comme dans notre article \cite{Brjuno}, les \'el\'ements de la th\'eorie classique des fractions continues qui nous seront utiles, et nous citons les estimations obtenues dans \cite{Brjuno} utiles à la démonstration du théorème \ref{t114}. Au \S\ref{t164} nous définissons la fonction de Wilton et décrivons son comportement local moyen en adaptant la démarche adoptée dans \cite{Brjuno} pour l'étude de la fonction de Brjuno.  Le \S\ref{t108} concerne la fonction $\fhi_1$ ; nous y déterminons notamment la relation liant $\fhi_1$ et la fonction $\W$. La méthode des équations fonctionnelles approchées, que nous avons adoptée à la suite de Wilton et qui constitue la source commune des deux parties de notre travail, est présentée sous une forme assez générale au \S\ref{par:eq-approchee}. Au \S\ref{derivabilite-phi2} nous étudions le comportement local de la primitive $\fhi_2$ de $2\fhi_1$ introduite par B\'aez-Duarte \textit{et al.} au \S4.2 de \cite{baez-duarte-all} ; nous montrons en particulier que ses points de dérivabilité sont précisément les points de convergence des séries $\fhi_1$ et $\W$ (dits \emph{points de Wilton}), et déterminons le comportement asymptotique de son module de continuité. Ces résultats sont enfin transférés à la fonction $A$ au \S\ref{t160}.

\emph{Deuxième partie : autour d'un théorème de Wilton, la Bretèche et Tenenbaum}. Au \S\ref{t120} nous présentons le contexte de ce théorème : il s'agit d'une classe de problèmes introduite par Davenport en 1937 (cf. \cite{davenport-1937-1}, \cite{davenport-1937-2}) pour élucider les relations entre convolution de Dirichlet et séries trigonométriques. Nous présentons également brièvement l'approche générale développée par la Bretèche et Tenenbaum dans \cite{breteche-tenenbaum}. Le \S\ref{par:psi1} présente une nouvelle démonstration de l'équation fonctionnelle approchée vérifiée par les sommes partielles de la série trigonométrique $\psi_1(x)$, découverte par Wilton. Notre démonstration est légèrement plus simple que celle de Wilton et fournit un meilleur terme d'erreur. Nous concluons cette partie au \S\ref{t149}.

\medskip

Nous n'avons pas rédigé de survol exhaustif des travaux liés au thème abordé dans cet article. Nous renvoyons cependant à d'autres approches de problèmes voisins par Jaffard (cf. \cite{jaffard}), Schoissengeier (cf. \cite{schoissengeier}), Aistleitner, Berkes et Seip (cf. \cite{aistleitner-berkes-seip}), Rivoal et Roques (cf. \cite{rivoal-roques}).

\section{Fractions continues}\label{par:fractions-continues}

\subsection{G\'en\'eralit\'es} 
Nous posons $X=]0,1[ \setminus\Rat $ et $\alpha (x)=\{1/x\}$ pour $x\in X$. Nous d\'efinissons les fonctions it\'er\'ees de $\alpha$ en posant $\alpha_0(x)=x$ et
pour $k\ge 1$, $\alpha_k(x)=\alpha \big(\alpha_{k-1}(x)\big)$.
Les fonctions $\alpha_k\,(k\ge0)$ sont d\'efinies sur $X$, \`a valeurs dans $X$. Lorsque $x$ est rationnel, les mêmes formules permettent de définir $\alpha_k(x)$ tant que $k\ioe K$, où $K$ est un entier, appelé \emph{profondeur} de $x$ (cf. \S\ref{profondeur_rationnel}) ; on a alors $\alpha_K(x)=0$.

\medskip
Rappelons le théorème de Ryll-Nardzewski (cf. \cite{MR0046583}, Theorem 2, p. 76) :  pour toute fonction $f\in L^1(0,1)$, on a $f\circ \alpha \in L^1(0,1)$ et 
\begin{equation}\label{eq:invariance-gauss}
\int_0^1 f\big(\alpha(t)\big) \frac{dt}{1+t} = \int_0^1 f(t) \frac{dt}{1+t}.  
\end{equation}

Pour $x\in X$, nous notons $a_0(x)=0$, et $a_k(x)= \lfloor 1/\alpha_{k-1}(x)
\rfloor \, (k \soe 1)$ les quotients incomplets de $x$ ; $p_k(x),
q_k(x) \, (k \soe 0)$ respectivement le num\'erateur et le d\'enominateur de la
fraction r\'eduite d'ordre $k$ de $x$. Avec la notation classique, on
a  
$$
[a_0(x);a_1(x),\dots,a_k(x)]=\frac{p_k(x)}{q_k(x)}.
$$
Les fonctions $a_k,p_k,q_k$ d\'efinies sur $X$  sont \`a valeurs dans $\Nat^*$ pour $k\ge 1$.  On a la relation 
\begin{equation}\label{eq:identite-alpha} 
\alpha_k(x)= -\frac{p_k(x)-xq_k(x)}{p_{k-1}(x)-xq_{k-1}(x)}.  
\end{equation}

On a, pour tout $k \soe 2$,
\begin{equation}\label{croissance_qk}
q_k=a_kq_{k-1}+q_{k-2} \soe q_{k-1}+ q_{k-2}\soe 2q_{k-2},
\end{equation}
et $q_k \soe F_{k+1}$  pour $k\soe 0$, o\`u $F_n$ est le $n$\up{e} nombre de Fibonacci ($F_0=0$, $F_1=1$ et $F_{n+1}=F_{n-1}+F_n$ pour $n \soe 1$).
On peut \'egalement (cf. \cite{Brjuno}, (13), p. 199) d\'eduire de \eqref{croissance_qk} la
majoration
\begin{equation}\label{eq:somme-qk}
q_0+q_1+\cdots +q_{k} \ioe 3q_{k},
\end{equation}
valable pour tout $k \in \Nat$.
Nous utiliserons également la majoration 
\begin{equation}\label{eq:somme-fibonacci}
 \sum_{k\ge 0} \frac{1}{F_{k+1}} \le 3,36.
\end{equation}

\subsection{Les fonctions $\beta_k$ et $\gamma_k$} 

Nous introduisons  les fonctions 
\begin{equation*}\label{t63}
 \beta_k(x)=\alpha_0(x)\alpha_1(x)\cdots\alpha_{k}(x)
\end{equation*}
(par convention, $\beta_{-1}=1$) et
 \begin{equation*}\label{defgamma}
 \gamma_k(x)=\beta_{k-1}(x)\log\big(1/\alpha_k(x)\big)\quad(x\in X, k\ge 0),
 \end{equation*}
de sorte que $\gamma_0(x)=\log 1/x$.
Rappelons les identit\'es
\begin{equation}\label{eq:identite-beta}
  \beta_k(x)=(-1)^{k-1}\bigl(p_k(x)-xq_k(x)\bigr)=\big|p_k(x)-xq_k(x) \big|=\frac{1}{q_{k+1}(x)+\alpha_{k+1}(x) q_k(x)} \quad(x\in X),
\end{equation}
les encadrements 
\begin{equation} \label{encabeta}
 \frac{1}{q_{k+1}+q_k} \ioe \beta_k \ioe \frac{1}{q_{k+1}} \quad (k \soe -1)
\end{equation}
et, (cf. \cite{Brjuno}, proposition 1)
\begin{equation}\label{enca_gamma}
-\frac{\log (2q_k)}{q_k}\ioe \gamma_k-\frac{\log q_{k+1}}{q_k}\ioe \frac{\log
2}{q_k}. 
\end{equation}
On dispose \'egalement de la majoration
 \begin{equation}\label{majobetakj}
 \beta_{i+j} \ioe \frac{1}{q_{i+1}F_{j+1}}\quad(i \soe -1, j \soe 0).
 \end{equation}

\subsection{Cellules}\label{par:cellules}
Soit $k\in\Nat$, $b_0=0$ et $b_1, \dots, b_k \in\Nat^*$. La cellule (de
profondeur $k$) $\cgot(b_1,\dots,b_k)$ est l'intervalle ouvert d'extr\'emit\'es
$[b_0;b_1,\dots,b_k]$ et $[b_0;b_1,\dots,b_{k-1},b_k+1]$.

Nous prolongeons par continuit\'e 
(sans changer de notation) les fonctions $a_k,p_k,q_k$ (et donc aussi $\alpha_k$) sur chaque cellule de profondeur $k$. 
Dans la cellule $\cgot(b_1,\dots,b_k)$, les fonctions $a_j$, $p_j$, $q_j$ sont ainsi constantes  pour $j \ioe k$ :
$$
a_j(x)=b_j, \quad \frac{p_j(x)}{q_j(x)}=[b_0;b_1,\dots,b_j] \quad (x \in \cgot(b_1,\dots,b_j)).
$$

Les extr\'emit\'es de $\cgot(b_1,\dots,b_k)$ sont
$$
\frac{p_k}{q_k} \quad \text{et} \quad \frac{p_k+p_{k-1}}{q_k+q_{k-1}}
$$
(dans cet ordre si $k$ est pair ; dans l'ordre oppos\'e si $k$ est impair).

\'Etant donn\'e un nombre irrationnel $x$ et un entier naturel $k$, il existe une unique cellule de profondeur $k$ qui contient $x$. Selon la proposition 4 de \cite{Brjuno}, la distance $\delta_k(x)$ de $x$ aux extr\'emit\'es de cette cellule satisfait \`a 
\begin{equation}\label{eq:majo-distance}
\delta_k(x)\le \frac{1}{q_k(x) q_{k+1}(x)}.
\end{equation}

\subsection{Dérivées des fonctions $\alpha_k$ et $\gamma_k$ dans une cellule de profondeur $k$}

Nous utiliserons les relations suivantes, valables dans une cellule de profondeur $k$.
\begin{align*}
\alpha'_k &=(-1)^k(q_k+\alpha_kq_{k-1})^2 \expli{cf. \cite{Brjuno}, (34), p. 207}\\
\gamma'_k &= (-1)^{k-1}q_{k-1}\log(1/\alpha_k)+(-1)^{k-1}/\beta_k \expli{cf. \cite{Brjuno}, (36), p. 208}.
\end{align*}

\subsection{Profondeur et \'epaisseur d'un segment}\label{profondeur} 

Soit $I=[a,b]$ un segment inclus dans $]0,1[$, de longueur $h=b-a>0$ et d'extr\'emit\'es $a,b$ irrationnelles. Il existe un unique entier naturel $K$ tel que $I$ soit inclus dans une cellule de profondeur $K$, mais dans aucune cellule de profondeur $K+1$. On dit alors que $I$ est de profondeur $K$. 
La longueur d'un segment $I$ tend uniform\'ement vers $0$ quand sa profondeur $K$ tend vers l'infini. Dans le cas o\`u $x \in X$ est fix\'e et $I=[x-h/2,x+h/2]$ avec $x\pm h/2\in X$, la profondeur $K=K(x,h)$ de $I$ tend vers l'infini quand $h$ tend vers $0$ (cf paragraphe 5.4 de \cite{Brjuno}). 

Nous d\'efinissons \'egalement l'\'epaisseur de $I$ comme le nombre de cellules de profondeur $K+1$ qui ont une intersection non vide avec $I$, o\`u $K$ est la profondeur de $I$. L'\'epaisseur de $I$ est donc un nombre entier sup\'erieur ou \'egal \`a $2$.

\subsection{Profondeur d'un nombre rationnel}\label{profondeur_rationnel}

Soit $r$ un nombre rationnel, $0<r<1$, mis sous forme irr\'eductible $p/q$. Il peut s'\'ecrire d'une et une seule fa{\c c}on sous la forme
$$
r=[0; b_1,\dots,b_k]
$$
avec $k \in \Nat^*$, $b_i \in \Nat^*$ pour $1\ioe i\ioe k$, et $b_k \soe 2$. Nous
dirons alors que $r$ est de profondeur $k$. On a alors $\alpha_k(r)=0$. Par
convention, $0$ et $1$ sont de profondeur $0$.  

\'Ecrivons $[0; b_1,\dots,b_{k-1}]$ sous forme r\'eduite $p_{k-1}/q_{k-1}$ (si $k=1$, on a $q_0=1$). Le nombre rationnel $r$ est une extr\'emit\'e de deux cellules de profondeur $k$ (qui sont donc contig\"ues) :
$$
\cgot \quad \text{d'extr\'emit\'es} \quad [0; b_1,\dots,b_{k-1},b_k-1]=\frac{p-p_{k-1}}{q-q_{k-1}}  \quad \text{et} \quad r \, ;
$$
$$
\cgot' \quad \text{d'extr\'emit\'es} \quad r \quad \text{et} \quad [0; b_1,\dots,,b_{k-1},b_k+1]=\frac{p+p_{k-1}}{q+q_{k-1}} .
$$

Dans le  paragraphe 5.3 de \cite{Brjuno}, nous avons \'etabli l'inclusion 
$$
]r-\frac{2}{3q^2},r+\frac{2}{3q^2}[ \, \setminus \, \{r\} \subset \cgot \cup \cgot'.
$$

Enfin, observons que les profondeurs des nombres rationnels appartenant à une cellule de profondeur $K$ sont supérieures à $K$.

\subsection{Comportement en moyenne de $\gamma_k$}
Nous restituons ici certaines majorations
de l'int\'egrale $\int_I\gamma_k(t)dt$ 
obtenues dans \cite{Brjuno} (propositions 7, 9 et 11).

\begin{prop}\label{cas_kinfK}
Soit $I$ un segment inclus dans $]0,1[$ de longueur $h$, avec $0<h\le e^{-2}$,
inclus dans une cellule de profondeur $K$ avec $K\in\Nat$. Pour $x\in I$,
$k\in\Nat$ tel que $k<K$, on a
$$
 \int_I \lvert \gamma_k(t)-\gamma_k(x)\rvert dt \ioe \frac{3}{2} q_{k+1} h^2,
$$
o\`u $q_{k+1}$ d\'esigne la valeur constante de cette fonction sur $I$. 
\end{prop}

\begin{prop}\label{cas_kegalK}
Soit $I$ un segment inclus dans $]0,1[$ de longueur $h$, avec $0<h\le e^{-2}$,
de profondeur $K$, dont le milieu $x$ et les extr\'emit\'es sont
irrationnels.  
On a
$$
\int_I \gamma_K(t) dt \ioe 8h \gamma_K(x)+\frac{h}{q_K} (6 \log q_K +4),
$$
o\`u $q_K$ d\'esigne la valeur constante de cette fonction sur $I$.
\end{prop}

\begin{prop}\label{cas_ksupK}
Soit $I$ un segment inclus dans $]0,1[$ de longueur $h$, avec $0<h\le e^{-2}$,
de profondeur $K$, d'\'epaisseur $E$, et dont le milieu $x$ et les extr\'emit\'es sont
irrationnels.
Pour $k \in \Nat$ tel que $k> K$ on a
$$
\int_I\gamma_{k}(t)dt \ioe \frac{h}{F_{k-K}}\Bigl (
\frac{72[E>2]}{q_K}+\frac{12[E=2]\log q_{K+2}(x)}{q_{K+1}(x)}+\frac{12[E=2]\log
q_{K+3}(x)}{q_{K+2}(x)}\Bigr ),
$$
o\`u $q_K$ d\'esigne la valeur constante de cette fonction sur $I$. 
\end{prop}

\section{La fonction de Wilton et son comportement local moyen}\label{t164}

\subsection{D\'efinition, rappels}\label{t94}

Dans \cite{Brjuno}, nous avons \'etudi\'e la fonction de Brjuno 
\begin{equation*}
 \Phi(x)=\sum_{k\ge 0} \gamma_k(x).
\end{equation*}
On a
\begin{align*}
 \int_0^1 \gamma_k(t)dt & \le \frac{1}{F_{k+1}} \int_0^1 \log\big(1/\alpha_k(t)\big)dt 
 \expli{d'après \eqref{encabeta} et \eqref{majobetakj}} \\
 & \le \frac{2}{F_{k+1}} \int_0^1 \log\big(1/\alpha_k(t)\big)\frac{dt}{1+t} \\ 
 & = \frac{2}{F_{k+1}} \int_0^1 \log(1/t)\frac{dt}{1+t} \expli{d'après \eqref{eq:invariance-gauss}} \\
 &\le  \frac{2}{F_{k+1}}.
\end{align*}
Ainsi la série à termes positifs $\Phi(x)$ est convergente dans $L^1(0,1)$, et donc presque partout. On appelle respectivement nombres de Brjuno et nombres de Cremer, ses points irrationnels de convergence, respectivement de divergence. 

\`A présent pour $x\in X=]0,1[\setminus \Rat$, nous posons : 
\begin{equation}\label{t112}
\Wcal(x)=\sum_{k\soe 0}(-1)^k\gamma_k(x).
\end{equation}

Si $x=r\in]0,1[$ est rationnel, les valeurs $\alpha_k(r)$ ne peuvent \^etre d\'efinies par la formule $\alpha_k(x) =\alpha\bigl (\alpha_{k-1}(x)\bigr)$ que pour $k\ioe K=K(r)$, o\`u $K(r)$ est la profondeur de $r$, définie au \S\ref{profondeur_rationnel}. On pose alors
\begin{equation}\label{t113}
\Wcal(r)=\sum_{0\ioe k <K}(-1)^k\gamma_k(r).
\end{equation}

Par exemple, $\W(1/k)=\log k$ si $k$ est entier, $k\soe 2$. Par convention, $\W(0)=0$. Enfin, on prolonge la d\'efinition de $\Wcal(x)$ \`a $x\in \Real$ par p\'eriodicit\'e : $\Wcal(x)=\Wcal(\{x\})$. 

\smallskip

La \emph{fonction de Wilton\footnote{En l'honneur du math\'ematicien australien John Raymond Wilton (1884-1944) qui a le premier introduit la s\'erie $\W(x)$ (cf. \cite{Wilton}, assertion $(20_{III})$, p. 221). L'article nécrologique \cite{MR0014998} donne un aper{\c c}u de sa vie et de son oeuvre.}} est la somme de la s\'erie \eqref{t112}, \eqref{t113} en tout point o\`u cette s\'erie converge (en particulier en tout point rationnel). Nous appelons \emph{nombres de Wilton}  (resp. \emph{nombres de Wilton-Cremer}) les points irrationnels de convergence (resp. de divergence).

\begin{prop}\label{prop:wilton-brjuno}
 Si $x\in X$ est un nombre de Brjuno  alors la s\'erie $\W(x)$ est convergente,
et de plus on a
\begin{equation}\label{eq:wilton-brjuno}
  |\W(x)| \le \Phi(x). 
\end{equation}
En particulier l'in\'egalit\'e \eqref{eq:wilton-brjuno} est valable presque partout. 
\end{prop}
\dem 
L'in\'egalit\'e \eqref{eq:wilton-brjuno} r\'esulte directement de l'in\'egalit\'e triangulaire. 
\fin 

\smallskip

Ainsi, tout nombre de Brjuno est un nombre de Wilton (et tout nombre de Wilton-Cremer est un nombre de Cremer). L'inégalité \eqref{eq:wilton-brjuno}
garantit également que la fonction $\W$ de Wilton appartient à $L^1(0,1)$.  
Bien que nous n'en ayons pas l'usage dans cette étude, signalons que la fonction de Wilton est à oscillation moyenne bornée\footnote{En conséquence, la fonction $\W$ appartient à $L^p(0,1)$ pour tout $p\ge 1$.} sur $[0,1]$, autrement dit que 
\[
\sup_{I} \frac{1}{|I|}\int_I\left| \W(x)- \frac{1}{|I|}\int_I \W(t)dt\right| dx<\infty,\]
où la borne supérieure est prise sur tous les intervalles $I$ inclus dans $[0,1]$.  
 Il est possible de le déduire des résultats obtenus dans \cite{Marmi_Moussa_Yoccoz}, où Marmi, Moussa et Yoccoz établissent que la fonction de Brjuno est à oscillation moyenne bornée. 

\begin{prop}\label{prop:critere-W}
 Soit $x\in X$. La s\'erie $\W(x)$ est convergente si, et seulement si 
la s\'erie
\begin{equation*}
  \sum_{k\ge 0} (-1)^k \frac{\log q_{k+1}(x)}{q_k(x)}
\end{equation*}
est convergente.
\end{prop}
\dem
C'est une cons\'equence directe de l'encadrement \eqref{enca_gamma}
et du fait que la s\'erie $$\sum_{k\ge 0} \frac{\log q_k}{q_k}$$ est convergente.
\fin

\subsection{\'Equation fonctionnelle de $\W$}\label{t116}

\begin{prop}\label{prop:equation-fonctionnelle-W}
Soit $x\in X$. Le r\'eel $x$ est un nombre de Wilton si, et seulement si $\alpha(x)$ est un nombre de Wilton. Dans ce cas, on a  
\begin{equation}\label{eq:equation-fonctionnelle-W} 
  \W(x)=\log(1/x)-x\W\bigl(\alpha(x)\bigr) \, , 
\end{equation}
 et $\alpha_k(x)$ est un nombre de Wilton pour tout entier naturel $k$. 
\end{prop}
\dem
Cela r\'esulte de l'identit\'e 
\begin{align*}
 \W(x)&= \log(1/x)+x\sum_{k\ge 1} (-1)^{k} \alpha_1(x)\cdots \alpha_{k-1}(x)
\log 1/\alpha_{k}(x) \\
 &= \log(1/x)+x\sum_{k\ge 1} (-1)^{k} \alpha_0\bigl(\alpha(x)\bigr)\cdots \alpha_{k-2}\bigl(\alpha(x)\bigr)
\log 1/\alpha_{k-1}\bigl(\alpha(x)\bigr)\\
&=\log(1/x)-x\W\bigl(\alpha(x)\bigr). \fine  
\end{align*}

\medskip 
On a plus g\'en\'eralement la proposition suivante. 
\begin{prop}\label{prop:equation-fonctionnelle-generale-W}
Pour tout nombre de Wilton $x$, $K\in \Nat$, on a 
\begin{equation}\label{eq:equation-fonctionnelle-generale-W}
\W(x)=\sum_{k<K}(-1)^k
 \gamma_k(x)+(-1)^K\beta_{K-1}(x)\W\bigl (\alpha_K(x)\bigr).
\end{equation}
De plus l'identit\'e \eqref{eq:equation-fonctionnelle-generale-W}
reste valable lorsque $x$ est un nombre rationnel de profondeur sup\'erieure ou \'egale \`a $K$. 
\end{prop}
\dem 
On peut l'\'etablir directement \`a partir de la d\'efinition de $\W(x)$ 
ou le d\'emontrer par r\'ecurrence \`a partir de la proposition 
\ref{prop:equation-fonctionnelle-W}. 
\fin

\subsection{Primitive de la fonction de Wilton} 

On introduit \`a pr\'esent la fonction absolument continue
\begin{equation*}
\intW(x)=\int_0^x \W(t)dt, \quad (0\ioe x\ioe 1).
\end{equation*}
Dans les paragraphes qui suivent, nous \'etudions le comportement local
de $\intW$. Nous emploierons \`a plusieurs reprises l'identit\'e 
\begin{equation*}
\int_I \sum_{k\ge 0}(-1)^k\gamma_k(t) dt = \sum_{k\ge 0}(-1)^k\int_I\gamma_k(t) dt, \quad(I \textrm{ intervalle inclus dans }[0,1]), 
\end{equation*}
dont la validit\'e est garantie par la majoration 
\begin{equation*}
 \int_0^1 \sum_{k\ge 0}\gamma_k(t) dt <\infty.
\end{equation*}

\subsection{Comportement de $\intW$ au voisinage d'un nombre rationnel} 
Les r\'esultats et d\'emonstrations de ce paragraphe sont analogues \`a ceux du paragraphe 7 
de \cite{Brjuno}, qui concernait la fonction $\Psi$, int\'egrale de la fonction de Brjuno. Dans un souci de lisibilit\'e, nous avons restitu\'e 
la plupart des d\'etails. 

\begin{lem}\label{majo-Wilton-moyenne}
 On a 
\begin{equation*}
\int_0^x \big| \W\bigl(\alpha(t)\bigr)\big| dt \ll x \quad(0\le x \le 1). 
\end{equation*}
\end{lem}
\dem 
On a 
\begin{align*}
\int_0^x \big|\W\bigl(\alpha(t)\bigr) \big| dt& \le \int_0^x \Phi\bigl(\alpha(t)\bigr)dt\expli{d'apr\`es 
la proposition \ref{prop:wilton-brjuno}} \\
&\ll x \expli{d'apr\`es le lemme 4 de \cite{Brjuno}}.\fine   
\end{align*}

\begin{lem}
Pour $0<x\le 1$, on a 
\begin{align}
 \intW(x)&=x\log(1/x) +x +O(x^2) \, ,\label{eq:comportement-zero}\\
 \intW(1)-\intW(1-x)&=-x\log(1/x)-x+O\big(x^2\log (2/x)\big)\, .\label{eq:comportement-un}
\end{align}
\end{lem}
\dem Pour $0<x\le 1$, 
\begin{align*}
\intW(x)=
\int_0^x \W(t)dt& =\int_0^x \Big (\log (1/t)-t\W\bigl(\alpha(t)\bigr)\Big) dt\expli{d'apr\`es la proposition
\ref{prop:equation-fonctionnelle-W}}\\
&= \int_0^x \log(1/t) dt +O\Big(x\int_0^x \big\lvert\W\bigl(\alpha(t)\bigr)\big\rvert dt \Big) \\
&= x\log(1/x) +x +O(x^2) \expli{d'apr\`es le lemme \ref{majo-Wilton-moyenne}}. 
\end{align*}

Pour établir \eqref{eq:comportement-un}, il suffit de traiter le cas $0<x<1/2$. Nous utiliserons notamment l'estimation
\begin{equation}\label{t169}
\int_0^x\Phi(u)du\ll x\log (1/x) \quad (0<x<1/2),
\end{equation}
qui résulte de l'estimation (46), p. 214 de notre article \cite{Brjuno}. On a
\begin{align*}
 \intW(1) -\intW(1-x)&=\int_{1-x}^1 \W(t) dt\\
&=\int_{1-x}^1 \Bigl (-t\W\bigl(\alpha(t)\bigr)+\log(1/t)\Bigr )\, dt \expli{d'apr\`es la proposition
\ref{prop:equation-fonctionnelle-W}}\\
&=\int_0^{x/(1-x)}\Bigl (-\W(u)+O\big(u\Phi(u)\big)+O(u)\Bigr )\, du \expli{changement de variables $u=(1-t)/t$}\\
&=-\frac{x}{1-x}\log \Big(\frac{1-x}{x}\Big)-
\frac{x}{1-x}+O\big(x^2\log(1/x)\big) \expli{d'apr\`es \eqref{eq:comportement-zero} et \eqref{t169}}\\
&=-x\log(1/x)-x+O\big(x^2\log(1/x)\big).\fine   
\end{align*}

Rappelons que la fonction $\W$ a été prolongée aux nombres rationnels par la formule \eqref{t113}.
\begin{prop}\label{prop:W-rationnel}
Pour $r=p/q\in ]0,1[$, $(p,q)\in \Nat^*\times \Nat^*$, $\mathrm{pgcd}(p,q)=1$,  $h\in \mathbb{R}$, 
$q^2|h|<2/3$, on a 
\begin{equation*}
 \intW(r+h)-\intW(r)
=\frac{1}{q}|h|\log(e/q^2|h|)
+h\W(r)+O\big(qh^2\log(1/q^2|h|)\big). 
\end{equation*}
\end{prop}
\dem 
Soit $K \in \Nat^*$ la profondeur du nombre rationnel $r$. Comme $0<|h|<2/3q^2$, on sait d'apr\`es le \S \ref{profondeur_rationnel} que l'intervalle ouvert d'extr\'emit\'es $r$ et $r+h$ est inclus dans l'une des deux cellules de profondeur $K$ :

$\bullet$ $\cgot$ d'extr\'emit\'es $\frac{p-p_{K-1}}{q-q_{K-1}}$ et $r$ ;

$\bullet$ $\cgot'$ d'extr\'emit\'es $r$ et $\frac{p+p_{K-1}}{q+q_{K-1}}$,

o\`u $p_{K-1}/q_{K-1}$ est la $(K-1)$\up{e} r\'eduite de $r$.

Les trois nombres $\frac{p-p_{K-1}}{q-q_{K-1}}$, $r$ et $\frac{p+p_{K-1}}{q+q_{K-1}}$ se succ\`edent dans cet ordre si $K$ est pair, dans l'ordre inverse si $K$ est impair. Par cons\'equent,
\begin{equation*}
  (-1)^Kh>0 \Leftrightarrow r+h\in \cgot'.
\end{equation*}

D'apr\`es la proposition \ref{prop:equation-fonctionnelle-generale-W} nous avons 
\begin{align*}
 \intW(r+h)-\intW(r)
&= \int_r^{r+h} \W(t) dt \\
&=\sum_{k<K} (-1)^k \int_r^{r+h} \gamma_k(t) dt 
+ (-1)^K \int_r^{r+h} 
\beta_{K-1}(t) \W\bigl (\alpha_K(t)\bigr )dt. 
\end{align*}
Tout d'abord,  
\begin{align*}
\Big|\sum_{k<K} (-1)^k \int_r^{r+h} \gamma_k(t) dt-h\W(r)\Big|
&=\Big|\sum_{k<K}(-1)^k \int_r^{r+h} 
\big( \gamma_k(t)-\gamma_k(r)\big) dt \Big|\\
&\le \frac{3}{2} 
h^2 \sum_{k<K} q_{k+1}\expli{d'apr\`es la proposition  \ref{cas_kinfK}}\\
& \le 5 q h^2 \expli{d'apr\`es la majoration \eqref{eq:somme-qk}}.\\ 
\end{align*} 

Il reste \`a montrer  
\begin{equation}\label{queue_rationnel}
(-1)^K \int_r^{r+h} \beta_{K-1}(t) \W\big(\alpha_K(t)\big) dt 
=\frac{|h|}{q}\log(1/q^2|h|)
+\frac{|h|}{q}+O\big(qh^2\log(1/q^2|h|)\big). 
\end{equation}
La fonction $\alpha_K$ est d\'efinie sur chacune des cellules $\cgot$ et $\cgot'$. Suivant que $r+h$ appartienne \`a $\cgot$ ou \`a $\cgot'$, nous prolongeons $\alpha_K$ (par continuit\'e \`a droite ou \`a gauche suivant la parit\'e de $K$) respectivement \`a $\cgot \cup \{r\}$ ou $\cgot' \cup \{r\}$, et dans l'int\'egrale du premier membre de \eqref{queue_rationnel} nous effectuons le changement de variables $u=\alpha_K(t)$. Distinguons les deux cas.

\noindent \textit{Premier cas} : $(-1)^Kh >0.$ 
Dans ce cas, le segment d'extr\'emit\'es $r$ et $r+h$ est inclus dans $\overline{\cgot'}$. Dans $\cgot'$, on a $q_K(t)=q$, $p_K(t)=p$ et, d'apr\`es \eqref{eq:identite-alpha}
\begin{equation*}
  u=\alpha_K(t)=\frac{q_Kt-p_K}{-q_{K-1}t
  +p_{K-1}}.
\end{equation*}
La fonction $\alpha_K$ est d\'erivable sur $\cgot'$ et 
$\alpha_K'=(-1)^K(q_K+\alpha_K q_{K-1})^2$, de sorte que  
\begin{equation*}
(-1)^K \int_r^{r+h} \beta_{K-1}(t) \W\big(\alpha_K(t)\big) dt
=(-1)^{2K} \int_{\alpha_K(r)}^{\alpha_K(r+h)}
  \W(u)\frac{du}{(q_K+q_{K-1}u)^3}.  
\end{equation*}
On a $\alpha_K(r)=0$ et
\begin{align*}
  \alpha_K(r+h)
  &=\frac{q^2h}{-q_{K-1}p-q_{K-1}qh+p_{K-1}q}\\
&=\frac{q^2|h|}{1-|h|qq_{K-1}},
\end{align*}
car $p_{K-1}q-pq_{K-1}=(-1)^K=\sgn \, h$. Nous posons donc
\begin{equation}\label{def_x'}
  x'=\frac{q^2|h|}{1-|h|qq_{K-1}}.
\end{equation}
Comme $q_{K-1}\le q/2$, on a
\begin{align*}
(-1)^K \int_r^{r+h} \beta_{K-1}(t) \W\big(\alpha_K(t)\big) dt 
&=\frac{1}{q^3}\int_0^{x'} \W(u) \frac{d u}{\Big(1+u\frac{q_{K-1}}{q}\Big)^3} \\  
&= \frac{x'\log(1/x')}{q^3} +\frac{x'}{q^3}+O\Big( \frac{x'^2}{q^3}\log(1/x')\Big), 
\end{align*}
o\`u la derni\`ere \'egalit\'e r\'esulte de 
\eqref{eq:comportement-zero} et \eqref{t169}.

\noindent \textit{Deuxi\`eme cas} : $(-1)^Kh <0.$
Dans ce cas, le segment d'extr\'emit\'es $r$ et $r+h$ est inclus dans $\overline{\cgot}$. Dans cette cellule, on a
\begin{equation*}
  u=\alpha_K(t)=\frac{(q-q_{K-1})t-p+p_{K-1}}{-q_{K-1}t
  +p_{K-1}},
\end{equation*}
Par suite, $\alpha_K(r)=1$ et
\begin{align*}
  \alpha_K(r+h)
  &= \frac{(q-q_{K-1})\frac{p}{q}-(p-p_{K-1})+h(q-q_{K-1})}{-q_{K-1}\frac{p}{q}+p_{K-1}
  -q_{K-1}h}\\
  &=1+\frac{hq^2}{(-1)^K-qq_{K-1}h}=1-x,
\end{align*}
avec cette fois
\begin{equation}\label{def_x}
  x=\frac{q^2|h|}{1+|h|qq_{K-1}}.
\end{equation}
Nous avons donc, d'apr\`es \eqref{eq:comportement-un} et \eqref{t169},
\begin{align*}
(-1)^K \int_r^{r+h} \beta_{K-1}(t) \W\big(\alpha_K(t)\big) dt
&= -\frac{1}{q^3}\int_{1-x}^1 \W(u)du + O\Big( \frac{x^2}{q^3}\log(1/x)\Big)\\
& = \frac{x\log(1/x)}{q^3} +\frac{x}{q^3}+O\Big( \frac{x^2}{q^3}\log(1/x)\Big). 
\end{align*}
  
Dans les deux cas on a donc pour $y=x$ ou $x'$, 
\begin{equation}\label{eq:rationnel-final} 
(-1)^K \int_r^{r+h} \beta_{K-1}(t) \W\big(\alpha_K(t)\big) dt 
=  \frac{y\log(1/y)}{q^3} +\frac{y}{q^3}+O\Big( \frac{y^2}{q^3}\log(1/y)\Big).
\end{equation}
Comme $q_{K-1}<q/2$ et $|h|<2/3q^2$, nous avons
\begin{equation*}
  \log(1/y)=\log(1/q^2|h|) +O(q^2h),
\end{equation*}
et
\begin{equation*}
  y=q^2|h|+O(h^2 q^4),
\end{equation*}
Ces estimations ins\'er\'ees dans \eqref{eq:rationnel-final} entra\^inent bien 
\eqref{queue_rationnel}.\fin

\subsection{Comportement de $\intW$ au voisinage d'un nombre de Wilton-Cremer} 

\begin{prop}\label{prop:Wilton-Cremer}
 Soit $x\in X$ un point de Wilton-Cremer. Alors la fonction 
$\intW$ n'est pas d\'erivable en $x$. 
\end{prop}

\dem 
\'Etant donn\'e un nombre $x$ de Wilton-Cremer, nous allons construire une suite
de nombres r\'eels $(h_K)_{K\soe 0}$ tendant vers $0$ quand $K\vers \infty$, et telle que le taux d'accroissement  
\begin{equation}\label{eq:taux-accroissement-W}
\frac{\intW(x+h_K)-\intW(x)}{h_K} 
\end{equation}
n'ait pas de limite quand $K\vers\infty$. 

Dans la suite pour tout entier naturel $k$, $a_k,q_k,p_k$ d\'esignent les valeurs des fonctions correspondantes au point $x$.  
Pour tout entier impair $K$, nous choisissons un nombre irrationnel $x_K$
appartenant \`a la cellule $\cgot(a_1,\ldots,a_K,a_{K+1}+2)$, tel que $(x+x_K)/2$
soit irrationnel. On a la suite d'in\'egalit\'es
\begin{align*}
 [a_0;a_1,\ldots,a_K,a_{K+1}]
<x&<[a_0;a_1,\ldots,a_K,a_{K+1}+1]&\\&<[a_0;a_1,\ldots,a_K,a_{K+1}+2]<x_K<[a_0;a_1,\ldots,a_K,a_{K+1}+3],
\end{align*}
 de sorte que le segment $I_K=[x,x_K]$ est inclus
dans la cellule $\cgot(a_1,\ldots,a_K)$, est de profondeur $K$, et d'\'epaisseur $3$. 
Notant $h_K=x_K-x$, nous avons
\begin{equation}\label{eq:majo-hK}
0< h_K \le \frac{p_K}{q_K}-x\le \frac{1}{q_Kq_{K+1}},  
\end{equation}
de sorte que $h_K$ tend vers $0$ lorsque $K$ tend vers l'infini. 
Enfin remarquons que la fonction $q_{K+1}$ satisfait \`a 
\begin{equation} \label{eq:enca-qk1}
q_{K+1} \le q_{K+1}(t) \le (a_{K+1}+2)q_K+q_{K-1}\le 3q_{K+1} \quad(t\in I_K).    \end{equation}

D'apr\`es la proposition \ref{cas_ksupK}, on a  
\begin{align}
\Big|\sum_{k>K}(-1)^k\int_{I_K} \gamma_k(t) dt\Big|& \le\sum_{k>K}\int_{I_K} \gamma_k(t) dt\notag\\ 
& \le \frac{72 h_K}{q_K}\sum_{k>K} \frac{1}{F_{k-K}} \notag\\ 
& \le \frac{242 h_K}{q_K}, \label{eq:Cremer1} 
\end{align}
o\`u la derni\`ere in\'egalit\'e provient de \eqref{eq:somme-fibonacci} et pourvu que $q_Kq_{K+1}\soe e^2$ (ce qui est vrai d\`es que $K\soe 3$).
De plus,  
\begin{align}
  \Big|\sum_{k<K}(-1)^k\int_{I_K} \big( \gamma_k(t)-\gamma_k(x)\big) dt\Big|& \le
\sum_{k<K}\int_{I_K} \big|\gamma_k(t)-\gamma_k(x)\big| dt\notag\\ 
& \le \frac{3 h_K^2}{2}\sum_{k<K} q_{k+1} \expli{d'apr\`es la proposition 
\ref{cas_kinfK}}\notag\\ 
& \le 5 h_K^2q_K \expli{d'apr\`es la majoration
\eqref{eq:somme-qk}}\notag\\
&\le \frac{5 h_K}{q_{K+1}} \expli{d'apr\`es la majoration \eqref{eq:majo-hK}}.\label{eq:Cremer2}
\end{align}
Par ailleurs,  la fonction $\gamma_K$ est d\'erivable sur $I_K$ et $\gamma_K'(t)=q_{K-1} \log\big(1/\alpha_K(t)\big) +1/\beta_K(t)$. On a donc, pour $t\in I_K$,   
\begin{align*}
|\gamma_K'(t)|&\le q_{K-1}\log\big(a_{K+1}(t)+1\big)+q_K+q_{K+1}(t)\expli{d'apr\`es \eqref{encabeta}}\\ 
&\le q_Ka_{K+1}(t)
+q_{K-1}+2q_{K+1}(t)\\
&=  3 q_{K+1}(t) \le 9 q_{K+1} \expli{d'apr\`es  \eqref{eq:enca-qk1}},   
\end{align*}
et ainsi 
\begin{align}
\Big| (-1)^K\int_{I_K} \bigl (\gamma_K(t)-\gamma_K(x)\bigr ) dt\Big|&\le  \int_{I_K} \big| \gamma_K(t)-\gamma_K(x)\big| dt\notag\\ 
&\le 9q_{K+1} \int_{I_K} (t-x) dt\notag\\ 
& = \frac 92 q_{K+1}h_K^2 \le 
\frac{5h_K}{q_K}. \label{eq:Cremer3}  
\end{align}
\`A pr\'esent nous avons 
\begin{align*}
\intW(x+h_K)-\intW(x)&=
\sum_{k\ge 0} (-1)^k \int_{I_K}
\gamma_k(t) dt \\
&= h_K\sum_{k\le K} (-1)^k \gamma_k(x) 
+ \sum_{k\le K} (-1)^k\int_{I_K} \big(\gamma_k(t)-\gamma_k(x)\big)dt 
\\&\quad + \sum_{k>K} (-1)^K \int_{I_K}
\gamma_k(t) dt,  
\end{align*}
et compte tenu de \eqref{eq:Cremer1}, \eqref{eq:Cremer2} et \eqref{eq:Cremer3}, nous obtenons finalement  
\begin{equation} \label{eq:final-Cremer}
 \frac{\intW(x+h_K)-\intW(x)}{h_K}=  
\sum_{k\le K}(-1)^k\gamma_k(x)
+O\Big(\frac{1}{q_K} \Big). 
\end{equation}

Le m\^eme raisonnement, mais en retenant cette fois les valeurs paires de $K$, fournit une  suite $h_K$, de nombres n\'egatifs, telles que \eqref{eq:final-Cremer} est v\'erifi\'ee. Comme $x$ est un point de Wilton-Cremer, la somme figurant au membre de droite 
de \eqref{eq:final-Cremer} diverge lorsque $K$ tend vers l'infini. 
Cela montre que le taux d'accroissement \eqref{eq:taux-accroissement-W}
diverge lorsque $h\to 0$. \fin 

\medskip

Observons que nous n'avons pas exclu la possibilit\'e que $\Upsilon$ soit d\'erivable \`a droite et \`a gauche en un point de Wilton-Cremer.

\subsection{Comportement de $\intW$ au voisinage d'un nombre de Wilton}  

\begin{prop}\label{prop:lebesgue-intW}
Tout nombre de Wilton est un point de Lebesgue de $\W$. En particulier la fonction $\intW$ est d\'erivable en tout nombre $x$ de Wilton et 
$\intW'(x)=\W(x)$.  
\end{prop}
\dem 
Soit $x$ un nombre de Wilton tel que $0<x<1$. Il s'agit d'\'etablir  
\begin{equation*}
 \int_{x-h/2}^{x+h/2}\big| \W(t)-\W(x) \big| 
dt =o(h) \quad(h\to 0). 
\end{equation*}
Comme $\Rat$ est dénombrable, on peut supposer que les nombres $x\pm h/2$ sont irrationnels. Posons alors $I=[x-h/2,x+h/2]$ et notons $K$ la profondeur de ce segment. De plus nous posons $q_j=q_j(x)$ pour tout entier naturel $j$.

Nous avons
\begin{equation*}
\begin{split}
\int_I \lvert\W(t)-\W(x) \rvert dt
&\le \sum_{k<K} \int_I \lvert \gamma_k(t)-\gamma_k(x) \rvert dt
+ \int_I \gamma_K(t) dt
 \\
&\quad+\sum_{k> K} \int_I \gamma_k(t) dt+h\Big|\sum_{k\soe K} (-1)^k\gamma_k(x)\Big|.
\end{split}
\end{equation*}

Majorons chacun des quatre termes du second membre. On a
\begin{align*}
\frac{1}{h}\sum_{k<K} \int_I \lvert \gamma_k(t)-\gamma_k(x) \rvert dt &\ioe \frac{3h}{2} \sum_{k<K} q_{k+1} \expli{d'apr\`es la proposition \ref{cas_kinfK}}\\
& \ioe \frac 92q_K h \quad \expli{d'apr\`es la majoration \eqref{eq:somme-qk}}\\
&\ioe \frac{9}{q_{K+1}} \expli{d'apr\`es la majoration  \eqref{eq:majo-distance}
puisque $h/2<\delta_K(x)$}.
\end{align*}
Ensuite, 
\begin{align*}
\frac{1}{h}  \int_I \gamma_K(t) dt &\ioe  8 \gamma_K(x)+\frac{6 \log q_K
+4}{q_K} \quad \text{\footnotesize (d'apr\`es la proposition \ref{cas_kegalK})}.
\end{align*}
Par ailleurs 
\begin{align*}
\frac{1}{h} \sum_{k> K} \int_I \gamma_k(t) dt &\ioe \Bigl ( \frac{72}{q_K}+12\frac{\log q_{K+2}}{q_{K+1}}+12\frac{\log q_{K+3}}{q_{K+2}}\Bigr ) \sum_{k> K}\frac{1}{F_{k-K}} \quad \expli{d'apr\`es la proposition \ref{cas_ksupK}}\\
& \ioe  \frac{242}{q_K}+41\frac{\log q_{K+2}}{q_{K+1}}+41\frac{\log q_{K+3}}{q_{K+2}} \expli{d'apr\`es la majoration \eqref{eq:somme-fibonacci}}.
\end{align*}

Puisque $x$ est un nombre de Wilton, on a
\begin{equation*}
\sum_{k\ge K}(-1)^k\gamma_k(x)=o(1) \quad  \textrm { et }\quad 
\gamma_K(x)=o(1) \quad(K\to \infty) 
\end{equation*}
et donc aussi, d'apr\`es l'encadrement \eqref{enca_gamma}, 
\begin{equation*} 
 \frac{\log q_{K+1}}{q_K}=o(1) \quad(K\to\infty). 
\end{equation*}
Finalement, on obtient 
\begin{equation}\label{majolongue}
\frac{1}{h}\int_I \lvert\W(t)-\W(x) \rvert dt= o(1) \quad(K\to\infty). 
\end{equation}
On obtient bien la conclusion souhait\'ee puisque $K$ tend vers l'infini
lorsque $h$ tend vers $0$ (cf. paragraphe \ref{profondeur}).

\subsection{Module de continuit\'e de $\intW$} 
Maintenant et dans la suite, nous notons 
\begin{equation*}
\omega(h,\fhi)=\sup \{|\fhi(x)-\fhi(y)|, 0\le x,y\le 1, |x-y|\le h\} \quad(h>0), 
\end{equation*}
le module de continuit\'e d'une fonction $\fhi:[0,1]\rightarrow \Real$. 

\begin{prop}\label{prop:mod-continuite-intW}
On a 
\begin{equation*}
 \omega(h,\intW)=h\log(1/h) +O(h) \quad(h>0). 
\end{equation*}
\end{prop}

\dem 

D'apr\`es \eqref{eq:comportement-zero}, nous avons 
\begin{equation*}
\omega(h,\intW) \ge h\log(1/h)+O(h).   
\end{equation*}
Par ailleurs, d'apr\`es le th\'eor\`eme 3 de \cite{Brjuno}, en notant 
$\Psi(x)=\int_0^x \Phi(t)dt$, on a 
\begin{equation*} 
 \omega(h,\Psi)= h\log(1/h)+O(h). 
\end{equation*}
Donc, pour $0\le x<y\le 1$, $|x-y|\le h$, 
\begin{equation*}
 \Big|\int_x^y \W(t) dt\Big| \le \int_x^y \big|\W(t)\big| dt
\le \int_x^y \Phi(t) dt \le \omega(h,\Psi) \le h\log(1/h)+O(h).     
\end{equation*}
Par suite, 
\begin{equation*}
 \omega(h,\intW)\le h\log(1/h)+O(h), 
\end{equation*}
ce qui ach\`eve la preuve. \fin

\section{La fonction $\fhi_1$ et son lien avec la fonction de Wilton}\label{t108}

Rappelons que la série $\fhi_1(t)$ est d\'efinie pour $t\in\Real$ par 
\[
\fhi_1(t)= \sum_{n\ge 1} \frac{B_1(nt)}{n}. 
\] 
 
Baez Duarte \textit{et al.} ont montr\'e (proposition 6 de \cite{baez-duarte-all}) que $\fhi_1$ d\'efinit une fonction  appartenant \`a $L^2(0,1)$, et donc \`a $L^1(0,1)$, en utilisant la repr\'esentation de $\fhi_1$ en s\'erie trigonom\'etrique.   On peut \'etablir que la série $\fhi_1$ converge dans $L^2(0,1)$ sans recourir \`a cette repr\'esentation en invoquant le crit\`ere de Cauchy, l'identit\'e
\[
 \int_0^1 B_1(mt)B_1(nt)dt=\frac{(m,n)^2}{12mn} 
 \quad(m,n\in\Nat^*)
\]
(démontrée par Landau de manière élémentaire : cf. \cite{50.0119.02}, p. 203) et le fait que la s\'erie double $\sum_{m,n\ge 1} \tfrac{(m,n)^2}{(mn)^2}$ est convergente.

Le but des prochains paragraphes est d'établir une relation entre les fonctions $\fhi_1$ et 
$\W$. Il nous faut au préalable introduire deux nouvelles fonctions auxiliaires.

\subsection{Les fonctions $F$ et $G$}\label{t163}

Posons pour $x >0$
\begin{equation}\label{eq:def-F} 
 F(x)=\frac{x+1}{2}A(1)-A(x)-\frac{x}{2} \log x.
\end{equation} 
 La fonction $F$ est continue sur $]0,\infty[$ et se prolonge par continuit\'e en
$0$ en posant $F(0)=A(1)/2$\footnote{On a $A(1)=\log 2\pi-\gamma$, o\`u $\gamma$
d\'esigne la constante d'Euler (cf. d\'emonstration de la proposition 8 de
\cite{baez-duarte-all}).}.

Pour $x \in X$, nous posons 
\begin{equation}\label{eq:serieF}
 G(x)= \sum_{j\ge 0}(-1)^j \beta_{j-1}(x)
  F\big(\alpha_j(x)\big)
\end{equation}
et nous prolongeons la d\'efinition de $G$ aux rationnels de $[0,1[$ en posant  
$$
G(r)= \sum_{j\ioe K}(-1)^j \beta_{j-1}(r)F\big(\alpha_j(r)\big),
$$
si $r$ est un nombre rationnel $\in [0,1[$ de profondeur $K$. En particulier, 
$G(0)=F(0)=A(1)/2$. Enfin nous prolongeons $G$ \`a $\Real$ par $1$-p\'eriodicit\'e. 

Les trois propositions suivantes r\'esument les principales propri\'et\'es de la fonction $G\, : \, [0,1[\vers \Real$.

\begin{prop}\label{prop:G1}
  La fonction $G$ est born\'ee sur $\Real$, continue en tout irrationnel, continue à droite en $0$, et on a $G(1-0)=-A(1)/2$.
\end{prop}
\dem
On peut limiter l'étude à l'intervalle $[0,1[$.

Pour $j\in \Nat$ et $x \in [0,1[$, nous posons 
$$
G_j(x)=
\begin{cases}
(-1)^j\beta_{j-1}(x)F\big(\alpha_j(x)\big)&\text{si $x \in X$ ou si $x$ est un rationnel de profondeur $\soe j$}\\
0&\text{si $x$ est un nombre rationnel de profondeur $<j$.}  
\end{cases}
$$

Nous avons donc
$$
G(x) =\sum_{j\soe 0}G_j(x) \quad (0\ioe x < 1).
$$

D\'esignons par $\|f\|_{\infty}$ la borne sup\'erieure de la valeur absolue d'une fonction $f$ sur $[0,1[$. Comme $\|G_j\|_{\infty} \ioe \|F\|_{\infty}/F_{j+1}$, la s\'erie est normalement convergente. En particulier, $G$ est born\'ee sur $[0,1]$.

D'autre part, si $x \in X$, toutes les fonctions $G_j$ sont continues en $x$. Il en r\'esulte que $G$ est continue en $x$.

Ensuite $G_0=F$ est continue sur $[0,1[$ et pour $j\soe 1$,
$$
|G_j(x)|\ioe \|F\|_{\infty}x,
$$
donc $G_j$ est continue en $0$.

Enfin, quand $x$ tend vers $1$,
\begin{align*}
  G_0(x)=F(x) &\vers \; \;\,F(1)=0,\\
G_1(x)=-x F\bigl (\alpha(x)\bigr) =-x F\bigl((1-x)/x\bigr)&\vers -F(0)=-A(1)/2\quad \text{et, pour $j\soe 2$,}\\
|G_j(x)|\ioe \beta_1(x)\|F\|_{\infty}=\|F\|_{\infty}(1-x)&\vers \;\;\;0,
\end{align*}
donc $G(1-0)=-A(1)/2$.\fin

\begin{prop}\label{prop:G2}
  La fonction $G$ v\'erifie l'\'equation fonctionnelle
  \begin{equation}
    \label{t98}
G(x)=F(x)-x G\bigl (\alpha(x)\bigr)\quad (x\in ]0,1[).    
  \end{equation}

Plus g\'en\'eralement, si $K\in \Nat$ et si $x\in [0,1[$ est irrationnel ou
rationnel de profondeur  sup\'erieure ou \'egale \`a $K$,
\begin{equation}\label{eq:equation-fonctionnelle-G-gene}
 G(x)=\sum_{j<K} (-1)^j \beta_{j-1}(x)F\big(\alpha_j(x)\big)
+(-1)^{K}\beta_{K-1}(x)G\big(\alpha_{K}(x)\big). 
\end{equation}
\end{prop}
\dem Soit $K\in \Nat$ et $x\in [0,1[$ irrationnel ou rationnel de profondeur $K'\soe K$. On a
$$
G(x) =\sum_{j<K} (-1)^j \beta_{j-1}(x)F\big(\alpha_j(x)\big)+
\begin{cases}
 \sum_{j\soe K} (-1)^j \beta_{j-1}(x)F\big(\alpha_j(x)\big) &\text{si $x \not\in \Rat$},\\
 \sum_{K\ioe j\ioe K'} (-1)^j \beta_{j-1}(x)F\big(\alpha_j(x)\big) &\text{si $x \in \Rat$}\\
\end{cases}
$$

Si $x \not\in \Rat$ on a
$$ 
\sum_{j\soe K} (-1)^j \beta_{j-1}(x)F\big(\alpha_j(x)\big) =(-1)^K \beta_{K-1}(x)G\big(\alpha_{K}(x)\big).
$$

Si $x \in \Rat$ on a
\begin{align*}
 \sum_{K\ioe j\ioe K'} (-1)^j \beta_{j-1}(x)F\big(\alpha_j(x)\big) &=\sum_{j'\ioe K-K'} (-1)^{K+j'} \beta_{j'-1+K}(x)F\big(\alpha_{j'+K}(x)\big)\\
&=(-1)^K \beta_{K-1}(x)\sum_{j'\ioe K-K'} (-1)^{j'} \beta_{j'-1}\bigl(\alpha_K(x)\bigr)F\Bigl(\alpha_{j'}\bigl(\alpha_K(x)\bigr)\Bigr )\\ 
&= (-1)^{K}\beta_{K-1}(x)G\big(\alpha_{K}(x)\big)
\end{align*}
puisque $\alpha_{K}(x)\in [0,1[$ est un nombre rationnel de profondeur $K'-K$.\fin 

\begin{prop}\label{t102}
  Si $r=p/q\in ]0,1[$ est un nombre rationnel de profondeur $K\soe 1$, \'ecrit sous forme irr\'eductible,  alors

$\bullet$ si $K$ est pair, $G$ est continue \`a droite en $r$ et $G(r-0)=G(r)-\frac{A(1)}{q}$ ;

$\bullet$ si $K$ est impair, $G$ est continue \`a gauche en $r$ et $G(r+0)=G(r)+\frac{A(1)}{q}$.
\end{prop}
\dem \'Ecrivons $r=\frac{p}{q}=[a_0;a_1,\ldots,a_K]$ avec $K\ge 1$ et $a_K\soe 2$ et
consid\'erons  $\cgot=\cgot(a_1,\ldots,a_K)$ et
$\cgot'=\cgot(a_1,\ldots,a_K-1,1)$, cellules contig\"ues de profondeurs
respectives  $K$ et  $K+1$, auxquelles $r$ est adh\'erent. Si $K$ est pair (resp. 
impair), alors la cellule
$\cgot$ se situe \`a droite (resp. \`a gauche) de $r$. 
Nous allons montrer que 
\begin{equation*} 
\lim_{ \substack{x\to r \\ x\in\cgot} } G(x)=G(r) 
\end{equation*}
tandis que 
\begin{equation*}
\lim_{ \substack{x\to r \\ x\in\cgot'} } G(x)=G(r)+(-1)^K\Big( \frac{F(1)}{q}-\frac{2F(0)}{q}\Big),
\end{equation*} 
ce qui montrera le r\'esultat souhait\'e. 
Remarquons tout d'abord que tout nombre rationnel appartenant \`a la r\'eunion 
$\cgot\cup \cgot'$ est de profondeur sup\'erieure ou \'egale \`a $K+1$. 
Pour $j< K$ les fonctions $\alpha_j$ et $\beta_j$ sont continues sur $\overline{\cgot\cup \cgot'}$.  
Par ailleurs notons que pour $x\in \overline{\cgot\cup \cgot'}\subset
\cgot(a_1,\ldots,a_{K-1})$, 
\begin{equation*}
 \beta_{K-1}(x) 
=|p_{K-1}-x q_{K-1}|=\frac{1}{q}+o(1) \quad(x\to r).
\end{equation*}

$\bullet$ Calcul de $\lim_{ \substack{x\to r \\ x\in\cgot} } G(x)$

Pour $x\in\cgot$, on a $p_K(x)=p$, $q_K(x)=q$, et par cons\'equent, d'apr\`es les formules  \eqref{eq:identite-alpha} et \eqref{eq:identite-beta}, 
\begin{align*}
\alpha_K(x)&=-\frac{p-x q}{p_{K-1}-x q_{K-1}}=o(1),\\  
\beta_K(x)&=|qx-p|=o(1) \quad (x\to r, \, x\in\cgot).
\end{align*}

D'apr\'es l'\'equation fonctionnelle \eqref{eq:equation-fonctionnelle-G-gene}, nous avons 
\begin{align*}
 G(x)=
\sum_{j\le K} (-1)^j \beta_{j-1}(x) F\big(\alpha_j(x)\big)
+(-1)^{K+1}|qx-p|G\big(\alpha_{K+1}(x)\big)\quad(x\in\cgot).
\end{align*}
Comme $G$ est born\'ee sur $[0,1[$ et comme $F$ est continue sur le m\^eme intervalle, il suit  
\begin{equation}\label{eq:limite-cgot}
\lim_{ \substack{x\to r \\ x\in\cgot} } G(x)
=\sum_{j< K} (-1)^j \beta_{j-1}(r) F\big(\alpha_j(r)\big)
+\frac{(-1)^K}{q}F(0)= 
G(r).  
\end{equation}

$\bullet$ Calcul de $\lim_{ \substack{x\to r \\ x\in\cgot'} } G(x)$

Tous les nombres rationnels de $\cgot'$ ont une profondeur sup\'erieure ou \'egale \`a 
$K+2$. 
Pour $x\in\cgot'$, on a $p_K(x)=p-p_{K-1}$, $q_K(x)=q-q_{K-1}$ de sorte que 
\begin{align*}
\alpha_K(x)&=-\frac{p-p_{K-1}-x(q-q_{K-1})}{p_{K-1}-x q_{K-1}}=1+o(1),\\
 \beta_K(x)&=\big|p-p_{K-1}-(q-q_{K-1})x\big|
=\frac{1}{q}+o(1) \quad( x\to r, \, x\in\cgot').
\end{align*}

De plus, pour $x\in\cgot'$, $p_{K+1}(x)=p$ et $q_{K+1}(x)=q$, et
\begin{align*} 
 \beta_{K+1}(x)&=|qx-p|,\\
\alpha_{K+1}(x)&=-\frac{p-x q}{p-p_{K-1}-x(q-q_{K-1})}=o(1) \quad(x\to r, \, x\in\cgot') . 
\end{align*}
D'apr\`es l'\'equation fonctionnelle, on a 
\begin{align*}
 G(x)&=\sum_{j< K} (-1)^j \beta_{j-1}(x) F\big(\alpha_j(x)\big)
+(-1)^K\Big(\beta_{K-1}(x)F\big(\alpha_K(x)\big)
-\beta_{K}(x)F\big(\alpha_{K+1}(x)\big)\Big) \\
&\quad 
+(-1)^{K+2}|p-qx|G\big(\alpha_{K+2}(x)\big) \quad(x\in \cgot').  
\end{align*}
Il suit 
\begin{equation*} 
\lim_{\substack{x\to r\\ x\in\cgot'}}  G(x)= 
 \sum_{j< K} (-1)^j \beta_{j-1}(r) F\big(\alpha_j(r)\big)+
(-1)^K\Big( \frac{F(1)}{q}-\frac{F(0)}{q}\Big)
=G(r)+(-1)^K\Big( \frac{F(1)}{q}-\frac{2F(0)}{q}\Big).   \fine
\end{equation*}

\subsection{Relation entre $\fhi_1$ et $\W$}\label{t161}

Introduisons la fonction $1$-p\'eriodique $\delta:\Real\to\Real$ d\'efinie par  
\begin{equation}\label{t122}
\delta(x)=
\begin{cases}
0 \textrm{ si } x\in X, \\[6pt]
\displaystyle{\frac{(-1)^{K+1} A(1)}{2q}} \textrm{ si } x=p/q \in[0,1[ , \; (p,q)=1, \; x \textrm{ de profondeur $K$}.
\end{cases} 
\end{equation}

L'objet de ce paragraphe est la démonstration de la proposition suivante, qui est la partie du théorème \ref{t150} concernant la série $\fhi_1$.

\begin{prop}\label{prop:lien-fhi1-W}
Les s\'eries $\fhi_1(x)$ et $\W(x)$ convergent pour les mêmes valeurs de $x\in[0,1]$ et l'identit\'e
\begin{equation}\label{eq:identite-G-X}
\fhi_1(x) = -\frac12 \W(x) + G(x) +\delta(x)
\end{equation} 
est valable en tout point de convergence. En particulier, on a $\fhi_1(x) = -\frac12 \W(x)+G(x)$ presque partout.
 \end{prop} 

\subsubsection{\'Equation fonctionnelle v\'erifi\'ee par les sommes partielles de $\fhi_1(x)$}\label{par:fhi1}

La première étape de la démonstration de la proposition \ref{prop:lien-fhi1-W} consiste à établir une équation fonctionnelle approchée vérifiée par les sommes partielles de $\fhi_1$.

\begin{prop}\label{prop:eq-approchee-fhi1}
Pour $x \in ]0,1]$, $v\in\Real$  et $xv\soe 1$, on a
$$
\sum_{1\ioe m\ioe v} \frac{ B_1(m x)}{m} +x \sum_{1\ioe n\ioe x
v} \frac{ B_1\bigl (n\alpha(x)\bigr)}{n}=F(x)-\demi\log
(1/x)+\eps(x,v) ,$$
avec 
\[
\eps(x,v)\ll \frac{1}{xv}. 
\] 
\end{prop}
\dem Nous reproduisons la formule figurant en haut de la page 225 de \cite{baez-duarte-all} :
\begin{multline}\label{t119}
\sum_{1\ioe m\ioe v} \frac{ B_1(m x)}{m} +x \sum_{1\ioe n\ioe x
v} \frac{ B_1 (n/x)}{n}=\frac{x}{2}\int_0^v\{t\}^2t^{-2}dt +
\frac{1}{2}\int_0^{x v}\{t\}^2t^{-2}dt- \int_0^v\{t\}\{x
t\}t^{-2}dt\\
+\frac{x -1}{2} \log (1/x)+\frac{x -1}{2}\int_v^{x v}\{t\}t^{-2}dt+\frac{1}{2 x v}(\{x v\}-x\{v\})^2 +\frac{x -1}{2x v}(\{x v\}-x\{v\}).
\end{multline}

Comme $B_1 (n/x)= B_1\bigl (n\alpha(x)\bigr)$, on reconna\^it bien la formule annonc\'ee, avec
\begin{multline*}
  \eps(x,v)=-\frac{x}{2}\int_v^{\infty}\{t\}^2t^{-2}dt-\demi\int_{x v}^{\infty}\{t\}^2t^{-2}dt+\int_v^{\infty}\{t\}\{x t\}t^{-2}dt\\
+\frac{x -1}{2}\int_v^{x v}\{t\}t^{-2}dt+\frac{1}{2 x v}(\{x v\}-x\{v\})^2 +\frac{x -1}{2 x v}(\{x v\}-x\{v\}).
\end{multline*}

Chacun des six termes composant $\eps(x,v)$ est $\ll 1/x v$ si $v\soe 1/x$ (ce qui entra\^ine $v\soe 1$).\fin

\medskip

Signalons que la source de l'identit\'e \eqref{t119} est la troisi\`eme d\'emonstration de la loi de r\'eciprocit\'e quadratique, propos\'ee par Gauss en 1808 (cf. \cite{gauss}, \S 5). Eisenstein en 1844 (cf. \cite{028.0828cj}) en donna une pr\'esentation g\'eom\'etrique particuli\`erement intuitive \`a l'aide d'un comptage de points \`a coordonn\'ees enti\`eres dans un rectangle. L'article d'Eisenstein fut traduit par Cayley  et publi\'e en 1857 dans le \textit{Quarterly Journal of pure and applied Mathematics}, \'edit\'e par Sylvester (cf. \cite{cayley}). Trois ans plus tard, Sylvester publia une note au Comptes Rendus \cite{sylvester} o\`u il g\'en\'eralisait l'argument d'Eisenstein \`a un rectangle $[0,v]\times [0,x v]$ o\`u $v$ et $x$ sont des r\'eels quelconques. La relation \eqref{t119} d\'ecoule de l'identit\'e de Sylvester par sommation partielle (cf. \cite{baez-duarte-all}, p. 223-225).

\subsubsection{Sur les solutions d'une \'equation fonctionnelle approch\'ee}\label{par:eq-approchee} 

Dans ce paragraphe, nous \'etablissons un r\'esultat g\'en\'eral sur des fonctions satisfaisant une \'equation approch\'ee apparentée à \eqref{t115}. 
  
On consid\`ere une fonction  $f:[0,1[\times [1,\infty[ \to \Real$ pour laquelle on suppose qu'il existe des nombres r\'eels $C\ge 1$, $a\ge 1$, $b>0$ tels que 
\begin{description}
\item[i)] 
pour tous $x\in]0,1]$, $v \in \Real$, tels que $x^a v\ge C$,  on a 
\begin{equation}\label{eq:fonctionnelle-approche}
f(x,v)+x f\big(\alpha(x),x^a v\big)=-\frac12 \log(1/x)+F(x)+ \eps(x,v),
\end{equation}
  avec 
\[
\eps(x,v) \ll (x^a v)^{-b}\, ; 
\] 
\item[ii)] pour tout $v\ge 1$, $f(0,v)=0.$
\item[iii)] on a uniform\'ement pour $x \in ]0,1[$ et $v\ge 1$
\[
f(x, v) \ll \log(2v) \, ; 
\] 
\end{description}

Pour tout $x\in [0,1]$ fix\'e, on se pose la question de l'existence de 
\[
\Fgot(x)=\lim_{v\to\infty} f(x,v).  
\] 
On a d'abord $\Fgot(0)=0$, d'après la condition ii). Ensuite, si $0<x\ioe 1$, la condition i) garantit que $\Fgot(x)$ existe si et seulement si $\Fgot\big(\alpha(x)\big)$ existe, et que dans ce cas, on a 
\begin{equation}\label{eq:wf-exacte}
\Fgot(x)+x\Fgot\big(\alpha(x)\big)=- \frac12 \log(1/x)+F(x).
\end{equation} 
Les deux sous-paragraphes qui suivent \'etablissent la proposition suivante. 
\begin{prop}\label{prop:meta-th} 
Soit $f:[0,1[\times [1,\infty[ \to \Real$ satisfaisant aux conditions i), ii) et iii) ci-dessus. La limite $\Fgot(x)$ existe si et seulement si la s\'erie $\W(x)$ est convergente et on a alors
\begin{equation*} 
\Fgot(x)=  - \frac12 \W(x)+G(x)+\delta(x). 
\end{equation*}  
\end{prop}

$\bullet$ \textbf{Le cas $x$ rationnel} 

\begin{prop} 
Soit $x=\frac{p}{q}\in[0,1[$ un rationnel de profondeur $K$ mis sous forme irr\'eductible. La limite $\Fgot(x)$ existe et de plus, 
\begin{equation}\label{eq:formule:rationnel}  
\Fgot(x)=- \frac12\W(x)+G(x) +(-1)^{K+1}\frac{A(1)}{2q},
\end{equation} 
\end{prop} 
\dem 
La convergence de $f(x,v)$ se d\'emontre par r\'ecurrence sur la profondeur $K$ de $x$. Si $K=0$, la convergence est directement garantie par la condition ii) et $\Fgot(0)=0$. 
Si $K\ge 1$,  $\alpha(x)$ est de profondeur $K-1$, et l'existence de $\Fgot(x)$ se d\'eduit directement de \eqref{eq:wf-exacte}. 
Pour \'etablir \eqref{eq:formule:rationnel}, on consid\`ere les trois relations
\begin{align*}
\W(x)&=\log(1/x)-x\W\big(\alpha(x)\big),\\ 
\Fgot(x)&= F(x)-\frac12 \log(1/x) -x\Fgot\big(\alpha(x)\big),\\
G(x)&= F(x)-xG\big(\alpha(x)\big),         
\end{align*}
v\'erifi\'ees par tous les rationnels $x \in ]0,1[$. En posant $h=\W+2\Fgot-2G $, on a 
\[h(0)=-2G(0)=-A(1),\] et
$$
h(x) =-x h\bigl(\alpha(x)\bigr) \quad (x \in ]0,1[ \,\cap\, \Rat).
$$
 Si $x=p/q$ est de profondeur $K\soe 1$, on a donc
\begin{align*}
  h(x) &=(-1)^j\beta_{j-1}(x)h\bigl( \alpha_j(x)\bigr)\quad (0\ioe j\ioe K)\expli{par r\'ecurrence sur $j$}\\
&=(-1)^K\beta_{K-1}(x) h(0)\\
&=(-1)^{K+1}A(1)/q.\fine
\end{align*}

$\bullet$ \textbf{Le cas $x$ irrationnel} 
Dans ce sous-paragraphe, nous notons abusivement $\beta_j(x)=\beta_j$ et $\alpha_j(x)=\alpha_j$ pour $j\ge 1$ afin de garantir une meilleure lisibilit\'e. 
 L'it\'eration de l'\'equation fonctionnelle approch\'ee i) fournit la proposition suivante.
\begin{prop}\label{t84}
Soit $x \in X$ et $K\in \Nat\,\cup\,\{-1\}$ tels que $\beta_K(x)^av\soe C$.
On a
\begin{equation}\label{t85}
\begin{aligned}
f(x,v)&+(-1)^K\beta_K f(\alpha_{K+1},\beta_K^a v)  
\\ &=-\frac12
\sum_{j=0}^K(-1)^j\beta_{j-1}\log(1/\alpha_j)
+\sum_{j=0}^K(-1)^j\beta_{j-1}F(\alpha_j) +\sum_{j=0}^K(-1)^j\beta_{j-1}
\eps\bigl(\alpha_j,\beta_{j-1}^a v).
\end{aligned}
\end{equation}
\end{prop}
\dem 
Le cas $K=-1$ est trivial. Le cas $K=0$ correspond \`a \eqref{eq:fonctionnelle-approche}.  Maintenant, si \eqref{t85} est vraie au rang $K$, et si de plus $v\soe C/\beta_{K+1}^a$, alors l'identit\'e \eqref{eq:fonctionnelle-approche} appliqu\'ee \`a $\alpha_{K+1}$ et $\beta_K^a v$ nous donne 
\begin{equation}\label{t86}
f\big(\alpha_{K+1}, \beta_K^a v\big) +\alpha_{K+1} 
f\big( \alpha_{K+2}, \beta_{K+1}^a v  \big)=-\frac12 \log(1/\alpha_{K+1})+F(\alpha_{K+1})+\eps\bigl(\alpha_{K+1},\beta_{K}^a v).  
\end{equation}
En ajoutant \eqref{t85} et $(-1)^{K+1}\beta_K\times$\eqref{t86}, on obtient bien la relation \eqref{t85} au rang $K+1$ au lieu de $K$.\fin

\begin{prop} 
Si $x\in X$, $f(x,v)$ converge lorsque $v\to\infty$ si et seulement si la s\'erie $\W(x)$ est convergente. On a alors  
\begin{equation}\label{t99} 
\lim_{v\to\infty} f(x,v) = \Fgot(x)= -\frac12 \W(x) + G(x). 
\end{equation} 
\end{prop} 
\dem
Supposons d'abord que la s\'erie $\W(x)$ est convergente. Si $v\soe 1$, d\'efinissons l'entier $K=K(x,v)\soe -1$ par l'encadrement
$$
\beta_{K}^a v \soe C>\beta_{K+1}^a v.
$$
On a
$$
q_{K+2}(x)+q_{K+1}(x)\soe \frac{1}{\beta_{K+1}}>(v/C)^{1/a} ,
$$
donc $K(x,v)\vers \infty$ quand $v\vers \infty$. Maintenant,
\begin{align}
(-1)^K\beta_K f(\alpha_{K+1},\beta_K^a v)  
&\ll \beta_K\log (2 \beta_{K}^a v)\label{t92}\\
&\ioe \beta_K\log  (2C/\alpha_{K+1}^a)\notag\\
&\ll \beta_K(x)+\gamma_K(x)\notag\\
&=o(1) \quad (K\vers \infty),\notag
\end{align}
puisque la s\'erie $\W(x)$ est convergente. 

Ensuite,
\begin{align}
 \sum_{j=0}^K(-1)^j\beta_{j-1}\eps\bigl(\alpha_j,\beta_{j-1}^a v)&\ll \sum_{j=0}^K\beta_{j-1}\bigl (\beta_{j}^a v\bigr)^{-b}\notag\\
&\le \frac{1}{v^b}\sum_{j=0}^K\frac{\beta_{j-1}}
{\beta_j^{ab}}. 
\end{align} 
La s\'erie $\sum_{j\ge 0} \beta_{j-1}$ est absolument convergente et la suite $\{\beta_{j}^{-ab} \}_{j\ge 0}$ tend vers l'infini en croissant. Un th\'eor\`eme classique de Kronecker (cf. \cite{18.0212.04}) garantit alors que 
\[
 \sum_{j=0}^K\frac{\beta_{j-1}}
{\beta_j^{ab}}=o(1/\beta_K^{ab}) \quad(K\to\infty).  
\] 
On a donc 
\begin{align} 
 \sum_{j=0}^K(-1)^j\beta_{j-1}\eps\bigl(\alpha_j,\beta_{j-1}^a v)
&=o(1/( \beta_K^a v)^b) \quad(K\to\infty)  \label{t91} \\
&=o(1) \quad (v \vers \infty).  \notag
\end{align} 
Par cons\'equent, la proposition \ref{t84} montre que $f(x,v)$ converge lorsque $v$ tend vers l'infini, et  a pour limite $-\frac12\W(x)+G(x)$.

\smallskip

Pour la r\'eciproque, on suppose que $f(x,v)$ converge lorsque $v\to\infty$. Si
$K\soe 1$ on pose $v=v(x,K)=C/\beta_{K}^a$, de sorte que
$v(x,K) \vers \infty$ quand $K\vers \infty$. Les estimations \eqref{t92}
et \eqref{t91} ci-dessus m\`enent \`a la conclusion que la s\'erie $\W(x)$
converge et a pour somme $-2\Fgot(x)+2G(x)$.\fin

\bigskip 

\subsubsection{Conclusion de la preuve de la proposition \ref{prop:lien-fhi1-W}} 
La somme partielle $\sum_{n\le v} B_1(nx)/n$ de la s\'erie $\fhi_1(x)$  satisfait aux points i), ii) et iii) \'enonc\'es dans le paragraphe \ref{par:eq-approchee}. 
En effet, le point ii) est trivial. En ce qui concerne le point iii), on a pour $x \in\Real$ et $v\ge 1$,  
\begin{equation*}
\sum_{n\le v} \frac{B_1(nx)}{n} \le \sum_{n\le v} \frac{1}{n} 
\ll \log(2v).  
\end{equation*} 
Enfin le point i) est v\'erifi\'e par la somme partielle de $\fhi_1(x)$ avec les param\`etres $a=b=1$ d'apr\`es la proposition \ref{prop:eq-approchee-fhi1}. La conclusion découle donc de la proposition \ref{prop:meta-th}.

\section{La fonction $\varphi_2$ et son comportement local}\label{derivabilite-phi2}

Soit $B_2$ la deuxi\`eme fonction de Bernoulli d\'efinie par 
\begin{equation*}
B_2(t)=\{ t \}^2-\{ t \}+\frac16 \quad(t\in\Real).  
\end{equation*}
Nous posons  pour $\lambda\in\Real$,
\begin{equation}\label{eq:def-phideux}
\fhi_2(\lambda)= \sum_{n\ge 1} \frac{B_2(n\lambda)}{n^2}, 
\end{equation}
qui d\'efinit une fonction $1$-p\'eriodique. 
De l'identit\'e classique $B_2(x)=\int_0^x 2B_1(t)dt +B_2(0)$, 
nous d\'eduisons 
\begin{align*}
\fhi_2(\lambda)&= \sum_{n\ge 1} \frac{1}{n^2} 
\int_0^{n\lambda} 2B_1(t)dt + B_2(0)\sum_{n\ge 1}\frac{1}{n^2} \\
&= 2 \sum_{n\ge 1} \frac{1}{n} \int_0^\lambda B_1(nu)du +
\frac{\zeta(2)}{6}\\
& =2  \int_0^\lambda  \sum_{n\ge 1} \frac{B_1(nu)}{n}du +
\frac{\zeta(2)}{6},
\end{align*}
l'interversion des signes $\sum$ et $\int$ dans la derni\`ere \'egalit\'e 
r\'esultant de la convergence dans $L^1(0,1)$ de la série $\fhi_1(u)$. Ainsi, pour $\lambda\in\Real$, 
\begin{equation}\label{formule-integrale-phideux}
\fhi_2(\lambda) 
=2 \int_0^\lambda \fhi_1(t)dt +\frac{\zeta(2)}{6},
\end{equation}
ce qui constitue la d\'efinition pour $\fhi_2$ adopt\'ee dans
\cite{baez-duarte-all}. 
En particulier, la fonction $\fhi_2$ est absolument continue. 
\begin{prop}\label{t100}
  On a
$$
\fhi_2(\lambda)-\fhi_2(0)=-\Upsilon(\lambda)+2\int_0^{\lambda}G(t)dt \quad (0\ioe \lambda \ioe 1).
$$
\end{prop}
\dem 

Cela r\'esulte imm\'ediatement de la formule \eqref{formule-integrale-phideux} et de la proposition \ref{prop:lien-fhi1-W}.\fin

\smallskip

Il est intéressant de retrouver partiellement le r\'esultat de la proposition 12 de \cite{baez-duarte-all}, obtenu alors par des moyens très différents (inversion de Mellin et équation fonctionnelle de la fonction d'Estermann).
 
\begin{prop}\label{prop:phi2-rat}
Soit $r=\frac{p}{q}$, $p,q\in \mathbb{N}^*$, $(p,q)=1$. La fonction $\varphi_2$ n'est pas d\'erivable 
en $r$ et 
\begin{equation}\label{eq:phi2-rat}
 \varphi_2\Big(\frac{p}{q}+h\Big)-\varphi_2\Big(\frac{p}{q}\Big)
= \frac{1}{q} |h|\log |h| +2\fhi_1(r)h+\bigl (2\log q-1+A(1)\bigr)|h|/q+o(h)\quad (h\to 0). 
\end{equation}
\end{prop} 
\dem 
Comme $\varphi_2$ est de p\'eriode $1$, on peut supposer que $r\in]0,1]$. Nous avons, d'apr\`es la proposition \ref{t100}, 
\begin{align*}
\varphi_2(r+h)-\varphi_2(r)
=\intW(r)-\intW(r+h)+2\int_r^{r+h} G(t) dt.  
\end{align*}

D'apr\`es la proposition \ref{prop:W-rationnel}, on a
\begin{equation}\label{t101}
 \intW(r)-\intW(r+h)
=-\frac{1}{q}|h|\log(1/|h|)
-h\W(r)+\frac{2\log q -1}{q}|h|+o(h) \quad (h\vers 0).
\end{equation}

D'apr\`es la proposition \ref{t102}, on a
\begin{equation}\label{t103}
\int_r^{r+h} G(t) dt=hG(r)+\frac{\sgn \, h+(-1)^{K+1}}{2}\cdot\frac{A(1)}{q}h+o(h) \quad (h\vers 0), 
\end{equation}
o\`u $K$ est la profondeur de $r$.

En ajoutant \eqref{t101} et le double de \eqref{t103}, on obtient
\begin{align*}
\varphi_2(r+h)-\varphi_2(r)&=-\frac 1q |h|\log(1/|h|)+h\bigl
(2G(r)-\Wcal(r)+(-1)^{K+1}A(1)/q\bigr)\\ &\quad+\bigl (2\log
q-1+A(1)\bigr)|h|/q+o(h), 
\end{align*}
ce qui donne bien le r\'esultat, compte tenu de \eqref{eq:identite-G-X}.\fin 

\begin{prop}\label{prop:phi2-wilton-cremer} 
 Soit $\lambda\in \Real$ un nombre de Wilton-Cremer. La fonction
$\varphi_2$ n'est pas d\'erivable en $\lambda$. 
\end{prop}
\dem 
Par $1$-p\'eriodicit\'e de $\varphi_2$, on peut supposer que $\lambda\in ]0,1[$.  
On a d'apr\`es la proposition \ref{t100}, 
\begin{equation*}
\frac{\varphi_2(\lambda+h)-\varphi_2(\lambda)}{h}
= - \frac{\intW(\lambda+h)-\intW(\lambda)}{h} 
+\frac{2}{h} \int_\lambda^{\lambda+h} G(t) dt.
\end{equation*}
Comme la fonction $G$ est continue en $\lambda$ (proposition 
\ref{prop:G1}), nous avons 
\begin{equation*}
\frac{\varphi_2(\lambda+h)-\varphi_2(\lambda)}{h}
= - \frac{\intW(\lambda+h)-\intW(\lambda)}{h} 
+2G(\lambda)+o(1) \quad(h\to 0).  
\end{equation*}
Or d'apr\`es la proposition \ref{prop:Wilton-Cremer}, le taux d'accroissement 
de la fonction $\intW$ diverge, ce qui fournit la conclusion attendue.  
\fin 

\begin{prop}\label{prop:phi2-wilton} 
 Soit $\lambda\in \Real$ tel que $\{\lambda\}$ est un nombre de Wilton. Alors $\lambda$ est un point de Lebesgue de $\fhi_1$. En particulier $\varphi_2$ est d\'erivable en $\lambda$ et 
$\varphi_2'(\lambda)=2\varphi_1(\lambda)$. 
\end{prop}
\dem 
Par $1$-p\'eriodicit\'e de $\fhi_1$ et $\varphi_2$, on peut supposer que $\lambda\in ]0,1[$. 
D'apr\`es la proposition \ref{prop:lien-fhi1-W} on a pour $h>0$
\begin{equation*}
\frac{1}{h} \int_{\lambda-h}^{\lambda+h}
\big|\fhi_1(t)-\fhi_1(\lambda)\big| dt 
\le \frac{1}{2}\int_{\lambda-h}^{\lambda+h}
\big|\W(t)-\W(\lambda)\big| dt+ \frac{1}{h}\int_{\lambda-h}^{\lambda+h}
\big|G(t)-G(\lambda) \big|dt.   
\end{equation*}
Comme la fonction $G$ est continue en $\lambda$, $\lambda$ est un point de Lebesgue de $G$. 
Par ailleurs, d'apr\`es la proposition \ref{prop:lebesgue-intW}, $\lambda$ est \'egalement un point de Lebesgue de $\W$. Par cons\'equent, 
\begin{equation*}
\frac{1}{h} \int_{\lambda-h}^{\lambda+h}
\big|\fhi_1(t)-\fhi_1(\lambda)\big| dt=o(1) \quad(h>0, h\to 0). \fine 
\end{equation*}

\smallskip

Nous pouvons d'autre part pr\'eciser la proposition 7 de \cite{baez-duarte-all}.

\begin{prop}\label{prop:phi2-module-continuite}
Le module de continuit\'e de $\varphi_2$ satisfait \`a 
\begin{equation*} 
\omega(h,\varphi_2) = h\log 1/h +O(h) \quad (0<h \ioe 1).  
\end{equation*}
\end{prop} 
\dem
Comme la fonction $G$ est born\'ee sur $[0,1]$, on a 
\begin{equation}
\Big| \int_x^y G(t) dt\Big|=O(h) \quad(0\le x<y\le 1, |x-y|\le h). 
\end{equation}
Puisque $2\varphi_1=-\W+2G$ presque partout, on en d\'eduit imm\'ediatement que 
\begin{equation*}
 \omega(h,\varphi_2)=\omega(h,\intW)+O(h) \quad(h>0). 
\end{equation*}
La conclusion r\'esulte alors directement de la proposition 
\ref{prop:mod-continuite-intW}. \fin 

\section{Comportement local de la fonction $A$ et fin de la preuve du théorème \ref{t114}}\label{t160}

\subsection{Relations entre $A$,  $\fhi_1$ et $\fhi_2$}\label{t162}

Dans ce paragraphe, nous donnons une démonstration élémentaire de la relation \eqref{t151} :
\begin{equation*}
 A(\lambda)=\demi\log \lambda
+\frac{A(1)+1}{2}-\lambda\int_{\lambda}^{\infty}\fhi_1(t)\frac{dt}{t^2}. 
\end{equation*}

Pour $\lambda>0$ et $N\in \Nat^*$, on a
\begin{align}
\int_{\lambda}^{\infty}\sum_{n\ioe N}\frac{B_1(nt)}{n}\,\frac{dt}{t^2} &=\sum_{n\ioe N}\int_{\lambda}^{\infty}B_1(nt)\frac{d(nt)}{(nt)^2}\notag\\
&=\sum_{n\ioe N}\int_{\lambda n}^{\infty}B_1(u)\frac{du}{u^2}\notag\\
&=\int_{\lambda}^{\infty}B_1(u)\min(N,\lfloor u/\lambda\rfloor)\frac{du}{u^2}\notag\\
&=\int_{\lambda}^{\lambda N}B_1(u)\lfloor u/\lambda\rfloor\frac{du}{u^2}+N\int_{\lambda N}^{\infty}B_1(u)\frac{du}{u^2}\label{t200}\\
&=\frac{1}{\lambda}\int_{\lambda}^{\lambda N}B_1(u)\frac{du}{u}-\int_{\lambda}^{\lambda N}B_1(u)\{ u/\lambda\}\frac{du}{u^2}+O(1/N\lambda^2)\notag\\
&=\frac{1}{\lambda}\int_{\lambda}^{\lambda N}B_1(u)\frac{du}{u}-\int_{\lambda}^{\lambda N}\{u\}\{ u/\lambda\}\frac{du}{u^2}+\demi\int_{\lambda}^{\lambda N}\{ u/\lambda\}\frac{du}{u^2}+O(1/N\lambda^2)\label{t153},
\end{align}
où l'on a utilisé la seconde formule de la moyenne pour majorer la deuxième intégrale de 
\eqref{t200}. 

Examinons successivement les trois intégrales figurant dans \eqref{t153}. On a d'abord
\begin{equation}\label{t154}
\int_{\lambda}^{\lambda N}B_1(u)\frac{du}{u}=\int_{\lambda}^{\infty}B_1(u)\frac{du}{u}+O(1/N\lambda).
\end{equation}
Ensuite,
\begin{align}
\int_{\lambda}^{\lambda N}\{u\}\{ u/\lambda\}\frac{du}{u^2} &=\int_{\lambda}^{\infty}\{u\}\{ u/\lambda\}\frac{du}{u^2}+O(1/N\lambda)\notag\\
&=A(1/\lambda)-\int_0^{\lambda}\{u\}\{ u/\lambda\}\frac{du}{u^2}+O(1/N\lambda),\label{t155}
\end{align}
où la dernière intégrale vaut
\begin{align}
\int_0^{\lambda}\{u\}\{ u/\lambda\}\frac{du}{u^2}&=  \frac{1}{\lambda}\int_0^{\lambda}\{u\}\frac{du}{u}\notag\\
&=1/\lambda +\frac{1}{\lambda}\int_1^{\lambda}\{u\}\frac{du}{u}\notag\\
&=1/\lambda +\frac{\log \lambda}{2\lambda}+\frac{1}{\lambda}\int_1^{\lambda}B_1(u)\frac{du}{u}.\label{t156}
\end{align}
Enfin,
\begin{align}
\int_{\lambda}^{\lambda N}\{ u/\lambda\}\frac{du}{u^2}&=\frac{1}{\lambda}\int_1^{N}\{ u\}\frac{du}{u^2}\notag\\
&=c_1/\lambda+O(1/N\lambda),\label{t157}
\end{align}
où $c_1=\int_1^{\infty}\{u\}/u^2du$.

En tenant compte de \eqref{t154}, \eqref{t155}, \eqref{t156} et \eqref{t157} dans \eqref{t153}, on obtient
\begin{equation}\label{t159}
\int_{\lambda}^{\infty}\sum_{n\ioe N}\frac{B_1(nt)}{n}\,\frac{dt}{t^2}=-A(1/\lambda)+\frac{\log \lambda}{2\lambda} +\frac{c_2}{\lambda} +O\big((\lambda^{-1}+\lambda^{-2})/N\big),
\end{equation}
où 
\begin{align}
c_2 &=1+\int_1^{\infty}B_1(u)\frac{du}{u}+\demi\int_1^{\infty}\{u\}\frac{du}{u^2}\label{t158}\\
&=1+\demi\int_1^{\infty}\{u\}^2\frac{du}{u^2}\expli{où l'on a intégré par parties la première intégrale de \eqref{t158}}\notag\\
&=\frac{A(1)+1}{2}.\notag
\end{align}

Comme la série de fonction périodique $\sum_n B_1(nt)/n$ converge dans $L^1(0,1)$, on peut faire tendre $N$ vers l'infini dans \eqref{t159} et obtenir ainsi
$$
\int_{\lambda}^{\infty}\fhi_1(t)\frac{dt}{t^2}=-A(1/\lambda)+\frac{\log \lambda}{2\lambda} +\frac{A(1)+1}{2\lambda},
$$
d'où l'identité \eqref{t151} en tenant compte de l'identité $A(\lambda)=\lambda A(1/\lambda)$.
 
\begin{prop}\label{prop:lien-A-phideux}
 On a pour tout $\lambda>0$, 
\begin{equation}\label{t78} 
 A(\lambda)=
\frac{1}{2} \log(\lambda) + \frac{1+A(1)}{2}
+ \frac{\fhi_2(\lambda)}{2\lambda}
-\lambda\int_\lambda^{\infty} \fhi_2(t) \frac{dt}{t^3}.
\end{equation}
\end{prop}
\dem

Cela résulte de \eqref{t151} par intégration par parties.\fin

\medskip 

On déduit immédiatement de la proposition \ref{prop:lien-A-phideux} et de l'identité  
$A(\lambda)=\lambda A(1/\lambda)$ les relations asymptotiques suivantes :
\begin{prop}
On a 
\begin{align}
 A(\lambda) &\sim \frac{1}{2} \log\lambda \quad(\lambda \to\infty), \label{t201}\\
 A(\lambda) & \sim  \frac{\lambda}{2} \log(1/\lambda)  \quad(\lambda \to 0)\label{t202}. 
\end{align}
\end{prop}

\subsection{Points de d\'erivabilit\'e} 

La proposition \ref{prop:lien-A-phideux} implique ainsi l'existence d'une fonction $f_1$ d\'erivable sur $]0,\infty[$ telle que 
\begin{equation}\label{t105} 
 A(\lambda)=\frac{\varphi_2(\lambda)}{2\lambda}+f_1(\lambda) \quad(\lambda>0).
\end{equation}
En particulier, $A$ et $\varphi_2$ ont les m\^emes points de d\'erivabilit\'e 
sur $]0,\infty[$. Compte tenu des r\'esultats de la section \ref{derivabilite-phi2}, nous obtenons le r\'esultat suivant.  
\begin{prop} 
Les points de d\'erivabilit\'e de $A$ sont les nombres de Wilton positifs.  
\end{prop}
Compte tenu de la proposition \ref{prop:critere-W}, nous avons bien \'etabli le th\'eor\`eme \ref{t114}. 

\subsection{Module de continuit\'e}

Dans ce paragraphe, nous donnons le comportement asymptotique du module de continuité de la fonction $A$.
\begin{prop}\label{t170}
 On a 
$$
\omega(A,h)=\frac 12h\log(1/h) +O(h)\quad (0<h\ioe 1).
$$
\end{prop}

\medskip

Les \'equations fonctionnelles $A(\lambda)=\lambda A(1/\lambda)$ et \eqref{t78} fournissent la relation
\begin{equation}\label{t79} 
 A(\lambda)=
\frac{\lambda}{2} \log(1/\lambda) + \frac{1+A(1)}{2}\lambda
+\frac{\lambda^2}{2}\varphi_2(1/\lambda)- \int_{1/\lambda}^{\infty}
\varphi_2(t) \frac{dt}{t^3}. 
\end{equation}

En particulier, on a
$$
A(h)-A(0)=\demi h\log(1/h) +O(h)\quad (0<h\ioe 1),
$$
ce qui prouve que $\omega(A,h)\soe \demi h\log(1/h) +O(h)$ pour $0<h\ioe 1$
tend vers $0$. Il nous reste \`a d\'emontrer l'in\'egalit\'e en sens inverse.

On a
$$
\omega(A,h)=\sup_{0<h'\ioe h}\max\bigl (\omega_0(A,h'),\omega_1(A,h')\bigr),
$$
avec
\begin{align*}
  \omega_0(A,h)&=\sup_{0< \lambda <1}|A(\lambda+h)-A(\lambda)|\\
  \omega_1(A,h)&=\;\sup_{\lambda \soe 1}|A(\lambda+h)-A(\lambda)|.
\end{align*}

\smallskip

Commen{\c c}ons par majorer $\omega_1(A,h)$. 
\begin{prop}
On a
$$
\omega_1(A,h)\ioe \demi h\log(1/h) +O(h) \quad (0<h\ioe 1).
$$
\end{prop}
\dem

On a, d'apr\`es \eqref{t78} et \eqref{t105},
$$
A(\lambda)=\frac{\varphi_2(\lambda)}{2\lambda} +f_1(\lambda),
$$
avec
\begin{align*}
  f_1(\lambda) &=\frac{1}{2} \log \lambda +\frac{1+A(1)}{2}-\lambda
\int_{\lambda}^{\infty}
\varphi_2(t) \frac{dt}{t^3},\\
f_1'(\lambda) &=\frac{1}{2\lambda}-\int_{\lambda}^{\infty}\varphi_2(t)
\frac{dt}{t^3}+\frac{\fhi_2(\lambda)}{\lambda^2}.
\end{align*}
En particulier, $f_1'(\lambda)=O(1)$ pour $\lambda \soe 1$, donc
$$
f_1(\lambda+h)-f_1(\lambda)=O(h) \quad (\lambda \soe 1,\, h>0).
$$

Maintenant,
$$
\frac{\varphi_2(\lambda +h)}{\lambda +h}-\frac{\varphi_2(\lambda)}{\lambda}=\frac{1}{\lambda +h}\bigl (\varphi_2(\lambda +h)-\varphi_2(\lambda)\bigr) -\frac{h}{\lambda(\lambda +h)}\fhi_2(\lambda).
$$
Le dernier terme est $O(h)$ pour $\lambda \soe 1,\, h>0$. D'autre part,
$$
|\varphi_2(\lambda +h)-\varphi_2(\lambda)|\ioe h\log(1/h) +O(h) \quad (0<h\ioe
1)
$$
d'apr\`es la proposition \ref{prop:phi2-module-continuite}. On a donc
\begin{equation}
  \label{t80}
  |A(\lambda +h)-A(\lambda)| \ioe \frac{1}{2(\lambda +h)}h\log(1/h) +O(h) \quad
(\lambda \soe 1,\, 0<h\ioe 1),
\end{equation}
et en particulier 
\begin{equation*}
\omega_1(A,h)\ioe \demi h\log(1/h) +O(h) \quad (0<h\ioe 1).\fine
\end{equation*}

\smallskip

Passons maintenant \`a l'\'etude de $\omega_0(A,h)$. D'apr\`es \eqref{t79}, on a
$$
A(\lambda) =\demi f_2(\lambda)+f_3(\lambda),
$$
avec
\begin{align*}
  f_2(\lambda) &=\lambda \log (1/\lambda)+\lambda^2\fhi_2(1/\lambda),\\
f_3(\lambda) &=\frac{1+A(1)}{2}\lambda- \int_{1/\lambda}^{\infty}\varphi_2(t)
\frac{dt}{t^3},\\
f'_3(\lambda) &= \frac{1+A(1)}{2}-\lambda\fhi_2(1/\lambda). 
\end{align*}

Comme $f'_3(\lambda)=O(1)$ pour $0<\lambda\ioe 2$, on a
$$
f_3(\lambda+h)-f_3(\lambda)=O(h) \quad (0<\lambda \ioe 1,\; 0<h\ioe 1).
$$

Pour démontrer la proposition \ref{t170}, il nous reste donc \`a montrer la proposition suivante.
\begin{prop}
On a
\begin{equation}
  \label{t81}
  |f_2(\lambda+h)-f_2(\lambda)|\ioe h\log(1/h) +O(h) \quad (0<\lambda \ioe 1,\;
0<h\ioe 1).
\end{equation}
\end{prop}
\dem

On a pour $0<\lambda \ioe 1$ et $0<h\ioe 1$, 
\begin{align*}
f_2(\lambda+h)-f_2(\lambda)&=\int_{\lambda}^{\lambda +h}\bigl (\log
(1/t)-1\bigr)dt+    \bigl ( (\lambda
+h)^2-\lambda^2\bigr)\fhi_2(1/\lambda)\\
&\quad+(\lambda +h)^2\Bigl (\fhi_2\bigl
(1/(\lambda +h)\bigr)-\fhi_2(1/\lambda)\Bigr)
\end{align*}
et
\[
\Big \lvert \int_{\lambda}^{\lambda +h}\bigl (\log (1/t)-1\bigr)dt+\bigl (
(\lambda +h)^2-\lambda^2\bigr)\fhi_2(1/\lambda)\Big \rvert \ioe h\min\bigl (\log
(1/\lambda),\log(1/h)\bigr) +O(h) .
\]

D'autre part,
\begin{align*}
  \Big \lvert\frac{1}{\lambda}-\frac{1}{\lambda +h}\Big \rvert&=\frac{h}{\lambda(\lambda +h)}\\
&\ioe \demi \quad (\lambda \soe \sqrt{2h}).
\end{align*}
Nous allons donc distinguer les deux cas $\lambda \soe \sqrt{2h}$ et $\lambda <\sqrt{2h}$.

Si $\lambda \soe \sqrt{2h}$, on a
\begin{align*}
(&\lambda +h)^2\Bigl \lvert\fhi_2\bigl (1/(\lambda
+h)\bigr)-\fhi_2(1/\lambda)\Bigr\rvert\\
 &\ioe (\lambda +h)^2\Bigl
(\frac{h}{\lambda(\lambda +h)}\log \Big(\frac{\lambda(\lambda +h)}{h}\Big)
+O\Bigl(\frac{h}{\lambda(\lambda +h)}\Bigr)\Bigr)\\
&=h\Bigl (1+\frac{h}{\lambda}\Bigr)\bigl (\log (1/h)+2\log \lambda +\log
(1+h/\lambda)\bigr)+
O\bigl (h (1+h/\lambda)\bigr)\\
&=h\log (1/h) +2h\log \lambda +O(h),  
\end{align*}
donc
$$
|f_2(\lambda+h)-f_2(\lambda)|\ioe h\log (1/h) +h\log \lambda +O(h)\ioe h\log (1/h) +O(h).
$$

Si $\lambda < \sqrt{2h}$, on a 
$$
(\lambda +h)^2\Bigl \lvert\fhi_2\bigl (1/(\lambda +h)\bigr)-\fhi_2(1/\lambda)\Bigr\rvert =O\bigl ((\lambda +h)^2\bigr)=O(h), 
$$
donc
on a encore
\begin{equation*}
|f_2(\lambda+h)-f_2(\lambda)|\ioe h\log (1/h) +O(h).\fine  
\end{equation*}

\section{Lien formel entre les s\'eries $\fhi_1(x)$ et $\psi_1(x)$}\label{t120}

Le développement en série de Fourier de la fonction $B_1$ est donn\'e par 
\[
B_1(t)= \frac{-1}{\pi} \sum_{m\ge 1} \frac{\sin(2\pi m t)}{m}, 
\] 
de sorte que l'on a formellement 
\begin{equation*}
\fhi_1(x)= -\sum_{m,n\soe 1}\frac{\sin (2\pi mn x)}{\pi mn} = \psi_1(x). 
\end{equation*}
 
Le probl\`eme de donner un sens analytique à cette identité formelle a \'et\'e plus g\'en\'eralement pos\'e par Davenport (\cite{davenport-1937-1}, \cite{davenport-1937-2}) en 1937 : \'etant donn\'e un couple de deux fonctions arithm\'etiques $(f,g)$ li\'ees par la relation 
\[f(n)=\sum_{d\mid n} g(d),\] 
pour quelles valeurs de $x$ l'identit\'e 
\[
\sum_{n\ge 1} \frac{g(n)}{n} B_1(nx)= -
\frac{1}{\pi} \sum_{m\ge 1} \frac{f(m)}{m}\sin(2\pi mx)
\]    
est-elle valide, au sens où les deux sommes sont convergentes et co\"incident? En g\'en\'eral, le problème le plus d\'elicat est celui de l'égalité des deux sommes ; il revient \`a montrer que 
\begin{equation}\label{eq:methode-davenport} 
\Delta(f,v,x):=\sum_{n\le v} \frac{g(n)}{n} B_1(nx) -
\frac{1}{\pi} \sum_{m\le v} \frac{f(m)}{m}\sin(2\pi m x)
\vers 0 \quad(v\to\infty).
\end{equation}   

Davenport traite notamment le cas du couple $(\delta,\mu)$, o\`u $\delta(n)=[n=1]$  et $\mu$ est la fonction de M\"obius. Sa m\'ethode  consiste \`a montrer la relation \eqref{eq:methode-davenport} en \'etudiant les variations et le comportement en moyenne de $x\mapsto \Delta(f,v,x)$. L'article \cite{breteche-tenenbaum} (pp. 4-5), ainsi que l'introduction de \cite{martin-2005} contiennent une description plus pr\'ecise de la m\'ethode de Davenport et de ses limites.  

\smallskip

 Dans \cite{breteche-tenenbaum} la Bret\`eche et Tenenbaum ont traité cette question grâce à la $P$-som\-ma\-tion. Rappelons ce dont il s'agit. 

Pour tout $n\in \Nat^*$, on note $P(n)$ le plus grand diviseur premier de $n$ (avec la convention $P(1)=1$). Une série de terme général $(u_n)_{n\soe 1}$ est dite $P$-sommable si

\smallskip

$\bullet$ pour tout $y\soe 1$ la série $\sigma(y)=\sum_{n\soe 1}u_n[P(n)\ioe y]$ converge ;

$\bullet$ la quantité $\sigma(y)$ a une limite quand $y$ tend vers l'infini.

\smallskip

Cette définition apparaît en 1957 au \S 2 de l'article \cite{MR0084629} de Duffin, mais l'étude approfondie de la $P$-sommation n'a véritablement commencé que 34 ans plus tard avec l'article \cite{MR1127146} de Fouvry et Tenenbaum. Cette étude s'appuie sur la connaissance des propriétés statistiques des entiers sans grand facteur premier, acquise depuis l'article fondateur de Dickman \cite{dickman} jusqu'à nos jours ; le chapitre III.5 de \cite{MR1366197} est un exposé pédagogique des résultats les plus importants de cette théorie.

L'approche de la Bretèche et Tenenbaum pour résoudre le problème de Davenport consiste à montrer au pr\'ealable la relation asymptotique 
\[
\nabla(f,v,x):= \sum_{P(n)\le v} \frac{g(n)}{n} B_1(nx) 
-\frac{1}{\pi} \sum_{P(m)\le v} \frac{f(m)}{m}\sin(2\pi m x)
=o(1) \quad(v\to\infty),
\] 
pour en d\'eduire \eqref{eq:methode-davenport}. Nous nous bornerons ici \`a renvoyer à
\cite{breteche-tenenbaum} et \cite{martin-2005} pour la description précise de cette m\'ethode, qui fournit de nombreux r\'esultats g\'en\'eraux et permet de traiter compl\`etement les cas embl\'ematiques des couples\footnote{Conform\'ement \`a l'usage, $\Lambda$ d\'esigne la fonction de Von Mangoldt, et $\one(n)=1$ pour tout $n\in\Nat^*$.} $( \log, \Lambda)$ et $(\tau,\one)$, ce dernier correspondant pr\'ecis\'ement \`a l'identit\'e \'etudi\'ee dans le pr\'esent travail. 

Notre m\'ethode consiste \`a montrer qu'en cas de convergence les sommes $\fhi_1(x)$ et $\psi_1(x)$ sont \'egales \`a une m\^eme troisi\`eme quantit\'e. Pour cela nous avons développé l'idée fondamentale de Wilton dans \cite{Wilton}, qui consistait à étudier la somme partielle $\psi_1(x,v)=-\sum_{n\ioe v}\tau(n)\sin (2\pi nx)/\pi n$ de la série $\psi_1(x)$ via une équation fonctionnelle approchée du type \eqref{eq:fonctionnelle-approche} :
\begin{equation}\label{t171}
 \psi_1(x,v)+x\psi_1\big(\alpha(x), x^c v \big)= -\frac{1}{2} \log(1/x)+F(x) +o(1) \quad(x\in]0,1], v\to\infty). 
\end{equation}
Nous avons vu au \S\ref{par:fhi1} que la somme partielle $\sum_{n\ioe v}B_1(nx)/n$ de la s\'erie $\fhi_1(x)$ v\'erifie une \'equation fonctionnelle approch\'ee de ce type avec $c=1$, o\`u $F$ est définie par \eqref{eq:def-F} ; cette équation fonctionnelle approchée était déjà présente implicitement dans \cite{baez-duarte-all}, et repose sur le même principe de symétrie que la démonstration géométrique de la loi de réciprocité quadratique. Au \S\ref{par:psi1}, nous précisons l'équation fonctionnelle de Wilton, en démontrant \eqref{t171} avec $c=2$ et une estimation améliorée du terme d'erreur $o(1)$. En outre, notre démonstration est sensiblement plus simple que celle de Wilton.

Ces équations fonctionnelles entra\^inent (cf. \S\ref{par:eq-approchee}) que le comportement des sommes partielles de $\fhi_1(x)$ et $\psi_1(x)$ suit asymptotiquement -- avec toutefois des termes correctifs aux points  rationnels-- celui de la solution $h(x)$ de l'\'equation fonctionnelle 
\begin{equation*}
h(x)+xh\big(\alpha(x)\big) =  -\frac{1}{2} \log(1/x)+F(x),
\end{equation*}
qui n'est autre que la fonction $-\frac12 \W(x) + \sum_{k\ge 0} (-1)^k \beta_{k-1}(x)F\big (\alpha_k(x)\big ).$  
\medskip

\section{\'Equation fonctionnelle v\'erifi\'ee par\\ les sommes partielles de $\psi_1(x)$}\label{par:psi1}

Rappelons la d\'efinition 
\[
\psi_1(x)=-\frac{1}{\pi}\sum_{n\ge 1} \frac{\tau(n)}{n}\sin(2\pi nx), 
\] 
o\`u $\tau(n)=\sum_{d\mid n}1$.

D\'eterminer une \'equation fonctionnelle approch\'ee pour les sommes partielles de $\psi_1(x)$ est bien plus ardu que pour celles de $\fhi_1(x)$. La raison en est que nous ne disposons pas d'une interprétation géométrique, analogue \`a celle donnant l'identité de Sylvester pour la fonction partie enti\`ere, qui nous permette  de relier simplement les sommes  
\[
\sum_{n\le v} \tau(n) \sin(2\pi nx) \textrm{ et } \sum_{n\le x^cv} \tau(n) \sin(2\pi n/x),
\] 
avec une valeur convenable de $c$.

La m\'ethode d\'evelopp\'ee par Wilton dans \cite{Wilton} repose sur l'utilisation de la formule sommatoire de Vorono\"i (cf. \cite{voronoi1}), employ\'ee par ce dernier pour donner une nouvelle démonstration de la majoration  
\begin{equation}\label{t127}
\Delta(x)=\sum_{n\le x} \tau(n)-x(\log x+2\gamma-1)  \ll x^{1/3+\eps}, 
\end{equation}
valable pour tout $\eps>0$, qu'il avait précédemment obtenue par une méthode élémentaire (cf. \cite{34.0231.03}). La formule de Voronoï, qui ramène l'étude d'une somme de la forme $\sum_n\tau(n)f(n)$ à une somme similaire où $f$ est remplacée par une certaine transformée intégrale de $f$, a permis à Wilton de découvrir (essentiellement) la proposition suivante.

\begin{prop}\label{prop:eq-approchee-psi1}
 On a pour $x\in]0,1]$, $v>0$, $x^2v\ge 2$, 
\begin{equation} \label{t2}
\frac{1}{\pi}\sum_{n\le  v} 
\frac{\tau(n)}{n} \sin(2\pi n x)
+\frac{x}{\pi}\sum_{n\le x^2v} \frac{\tau(n)}{ n} \sin(2\pi n /x)
 =-\frac{1}{2}\log x- F(x) +O\big((x^2v)^{-1/2}\log^2(x^2v)\big),
\end{equation} 
où $F$ est définie par \eqref{eq:def-F}.
\end{prop} 

La proposition \ref{prop:eq-approchee-psi1} découle du point (2.2) du Theorem 2, p. 223 de \cite{Wilton} (il faut prendre la partie imaginaire de la relation donnée par Wilton), à ceci près que le second membre de \eqref{t2} y apparaît sous la forme
$$
-\demi \log x + \Fgot (x) +O\big ((x^2v)^{-1/5}\big),
$$
où $\Fgot$ est une fonction continue sur $[0,1]$, non explicitée. 

Cela étant, plutôt que d'exposer les détails de la démarche de Wilton, nous avons choisi une autre voie, qui consiste au fond à utiliser la formule sommatoire de Voronoï sous la forme symétrique donnée par Nasim (cf. \cite{MR0284410}, Main theorem, p. 36). Cependant nous explicitons les détails nécessaires, sans supposer une étude préalable de \cite{MR0284410}.

Notre approche repose sur : la transformation de Mellin
\begin{equation}\label{t123}
Mf(s)=\int_0^{\infty}f(t)t^{s-1}dt \, ,
\end{equation}
la transformation de Fourier en cosinus
\begin{equation}\label{t124}
\Cgot f(x)=2\int_0^{\infty}f(t) \cos (2\pi tx)dt \, ,
\end{equation}
et quelques propriétés des fonctions cosinus intégral généralisées
\begin{equation}\label{t125}
\ci(a,z)=\int_z^{\infty}\cos t \cdot t^{a-1}dt\, ; \quad \Ci(a,z)=\int_0^{z}\cos t \cdot t^{a-1}dt
\end{equation}
(pour ces dernières, nous avons adopté les notations de \cite{MR2723248}, \S 8.21, p. 188).

Les trois sous-paragraphes suivants présentent les propriétés de \eqref{t123}, \eqref{t124} et \eqref{t125} que nous utilisons au \S\ref{t126} pour démontrer \eqref{t2}. 

\subsection{Rappels sur la transformation de Mellin}

Nous recommandons la lecture de l'appendice A (p. 231) de \cite{baez-duarte-all} et nous en rappelons ci-dessous les éléments qui nous seront utiles.

En désignant par $s$ la variable complexe, on note $\sigma=\Re s$ et $\tau=\Im s$. Si $-\infty\ioe a<b\ioe\infty$, on note\footnote{Nous conservons ici cette notation, aucune confusion avec la fonction de Wilton n'étant à craindre.}  $\Wcal(a,b)$ l'ensemble des fonctions complexes $f$ mesurables sur $]0,\infty[$ telles que
$$
\int_0^{\infty}|f(t)|t^{\sigma-1}dt<\infty \quad (a<\sigma<b).
$$
Si $f\in \Wcal(a,b)$, la transformée de Mellin $Mf$ définie par \eqref{t123} est holomorphe dans la bande $a<\sigma<b$. 

Nous considérons maintenant les transformées de Mellin qui interviennent dans notre argumentation, et rappelons, au \S\ref{t135}, la forme que prend la théorie de Plancherel dans le contexte de la transformation de Mellin.

\subsubsection{La fonction $A$}

On sait que $A\in \Wcal(-1,0)$ et que
$$
MA(s)=-\frac{\zeta(-s)\zeta(s+1)}{s(s+1)} \quad (-1<\sigma <0)
$$
(cf. \cite{baez-duarte-all}, proposition 10).

On en déduit que la fonction $A_1$ définie par $A_1(t)=A(t)/t$ ($=A(1/t)$) appartient à $\Wcal(0,1)$ et que
$$
MA_1(s)=\frac{\zeta(1-s)\zeta(s)}{s(1-s)} \quad (0<\sigma <1).
$$

Comme cette fonction appartient\footnote{Rappelons que pour $0<\sigma < 1$, $\tau\ge \tau_0$, 
$ \zeta(\sigma+i\tau) \ll \tau^{(1-\sigma)/2} \log(\tau)$ (cf. \cite{ivic} theorem 1.9 p. 25 par exemple).  
} 
 à $L^1(\sigma+i\Real,d\tau/2\pi)$ pour tout $\sigma\in]0,1[$, la formule d'inversion de Mellin nous donne 
\begin{equation}\label{t145}
A(x)=\int_{\Re s=\demi}\frac{\zeta(1-s)\zeta(s)}{s(1-s)}x^s\frac{d\tau}{2\pi},
\end{equation}
où l'égalité, \emph{a priori} valable presque partout, est vraie pour tout $x>0$ par continuité.

\subsubsection{Le reste dans le problème des diviseurs de Dirichlet}

Nous utiliserons également l'expression de la transformée de Mellin du reste $\Delta$, défini par \eqref{t127}, dans le problème des diviseurs de Dirichlet. \`A cette occasion, donnons l'énoncé général d'un principe classique, qui est une réciproque de la proposition 14 de \cite{baez-duarte-all}.
\begin{prop}\label{t129}
Soit $a<b\ioe c<d$, $f\in \Wcal(a,b)$, $g\in\Wcal(c,d)$. On suppose que $g-f$ est un polynôme généralisé
$$
P(t)=\sum_{\rho,k} c_{\rho,k}t^{-\rho}\log^kt,
$$
où la somme est finie, les $\rho$ sont des nombres complexes vérifiant $b\ioe \Re \rho\ioe c$, les $k$ sont des entiers naturels, et les $c_{\rho,k}$ des coefficients complexes.
Alors $Mf$ et $Mg$ sont les restrictions aux bandes $a<\sigma<b$ et $c<\sigma<d$, respectivement, d'une même fonction méromorphe dans la bande $a<\sigma<d$, dont la somme des parties polaires est
$$
\sum_{\rho,k} c_{\rho,k}\frac{(-1)^kk!}{(s-\rho)^{k+1}}.
$$
\end{prop}
\dem

En résumé, on va montrer que les fonctions $Mf$ et $Mg$ sont des prolongements méromorphes l'une de l'autre.

On a
\begin{equation}\label{t128}
f(t)+P(t)[t\soe 1]=g(t)-P(t)[t< 1].
\end{equation}

En notant $h(t)$ la valeur commune des deux membres de \eqref{t128}, on voit que $h\in \Wcal(a,b)$ (premier membre) et $h\in\Wcal(c,d)$ (second membre). On a donc $h\in \Wcal(a,d)$ ; la transformée de Mellin $Mh$ est holomorphe dans la bande $a<\sigma<d$, et
\begin{align*}
Mf(s) &= Mh(s) +\sum_{\rho,k} c_{\rho,k}\frac{(-1)^kk!}{(s-\rho)^{k+1}} \quad (a<\sigma<b),\\
Mg(s) &= Mh(s) +\sum_{\rho,k} c_{\rho,k}\frac{(-1)^kk!}{(s-\rho)^{k+1}} \quad (c<\sigma<d).\fine
\end{align*}

La proposition \ref{t128} s'applique à
\begin{align*}
f(t) &=\sum_{n\ioe t}\tau(n),\\
g(t) &=\Delta(t),\\
P(t)&=-t(\log t +2\gamma-1).
\end{align*}
On a $f\in\Wcal(-\infty,-1)$, $g\in\Wcal(-1,-1/3)$ (en utilisant l'estimation \eqref{t127}). Comme
\begin{align*}
Mf(s)&=\int_0^{\infty}\big (\sum_{n\ioe t}\tau(n)\big)t^{s-1}dt\\
&=\sum_{n\soe 1}\tau(n)\int_n^{\infty}t^{s-1}dt\\
&=\frac{\zeta^2(-s)}{-s}\quad(\sigma<-1),
\end{align*}
on en déduit que 
$$
M\Delta(s)=\frac{\zeta^2(-s)}{-s}\quad(-1<\sigma<-1/3).
$$

Par conséquent la fonction $\Delta_1$ définie par $\Delta_1(t)=\Delta(t)/t$ appartient à $\Wcal(0,2/3)$ et 
$$
M\Delta_1(s)=\frac{\zeta^2(1-s)}{1-s}\quad(0<\sigma<2/3).
$$

\subsubsection{La fonction $t\mapsto (\cos t )[t\ioe v]$}

Pour tout $v>0$, la fonction $t\mapsto (\cos t )[t\ioe v]$ appartient à $\Wcal(0,\infty)$ et sa transformée de Mellin est $\Ci(s,v)$. Par conséquent la transformée de Mellin de $t\mapsto (\cos 2\pi t x)[t\ioe v]$ est $(2\pi x)^{-s}\Ci(s,2\pi xv)$.

\subsubsection{La transformation de Mellin-Plancherel sur $L^2(0,\infty)$}\label{t135}

Si $f\in L^2(0,\infty)$, la formule \eqref{t123}, où l'intégrale doit être comprise comme $\lim_{T\vers \infty}\int_{1/T}^T$ dans $L^2(\demi+i\Real,d\tau/2\pi)$ définit un élément de cet espace, et l'application ainsi définie, dite transformation de Mellin-Plancherel, est une isométrie bijective entre espaces de Hilbert (cf. par exemple \cite{titchmarsh-fourier} \S 3.17 pp. 94-95). En particulier, si $f$ et $g$ appartiennent à $L^2(0,\infty)$, on a
$$
\int_0^{\infty}f(t)\overline{g(t)}\,dt =\int_{\sigma=1/2}Mf(s)\overline{Mg(s)}\,\frac{d\tau}{2\pi}
$$
(théorème de Plancherel).

Appliquons  le théorème de Plancherel à l'intégrale
\begin{equation}\label{t136}
I(x,v)=2\int_0^v \frac{\Delta(t)}{t}\cos  (2\pi tx)dt\, ,
\end{equation}
qui va jouer un rôle essentiel au \S\ref{t126}. On a
\begin{align}
I(x,v)&=2\int_0^{\infty}\Delta_1(t)\cos (2\pi tx)[t\ioe v]dt\notag\\
&=2\int_{\Re s=\demi}\frac{\zeta^2(1-s)}{1-s}\overline{\Ci(s,2\pi xv)}(2\pi x)^{s-1}\,\frac{d\tau}{2\pi}\notag\\
&=\int_{\Re s=\demi}\frac{\zeta^2(1-s)}{1-s}\Ci(1-s,2\pi xv)(2\pi x)^{s-1}\,\frac{d\tau}{\pi}.\label{t5}
\end{align}

\subsection{Rappels sur la transformation de Fourier en cosinus}

Si $f\in L^2(0,\infty)$, la formule \eqref{t124}, où l'intégrale doit être comprise comme $\lim_{T\vers \infty}\int_{0}^T$ dans $L^2(0,\infty)$, définit un élément de cet espace, et l'application $\Cgot$ ainsi définie est un opérateur unitaire involutif de $L^2(0,\infty)$, dite transformation de Fourier en cosinus.

La démonstration proposée au \S\ref{t126} est sous-tendue par le lien existant entre $\Cgot$ et la transformation de Mellin-Plancherel $M$ sur $L^2(0,\infty)$, lien qui est l'une des nombreuses façons de concevoir l'équation fonctionnelle de la fonction $\zeta$ de Riemann : 
\begin{quote}
si $f\in L^2(0,\infty)$, alors 
\begin{equation}\label{t130}
M(\Cgot f)(s) =\frac{\zeta(1-s)}{\zeta(s)}Mf(1-s),
\end{equation}
\end{quote}
égalité entre deux éléments de $L^2(\demi+i\Real,d\tau/2\pi)$ (\cad vraie pour presque tout $\tau$, et $\sigma=\demi$). En effet les deux membres de \eqref{t130} définissent des isométries de $L^2(0,\infty)$ dans $L^2(\demi+i\Real,d\tau/2\pi)$ et coïncident sur les fonctions $t\mapsto [0\le t\le a]$, dont les combinaisons linéaires constituent une partie dense de $L^2(0,\infty)$.

Comme application de \eqref{t130}, donnons deux relations entre la fonction d'au\-to\-cor\-ré\-la\-tion $A$ et le reste dans le problème des diviseurs de Dirichlet.
\begin{prop}\label{t167}
Pour $x > 0$, on a
\begin{align}
A(x) &=\int_0^x \Cgot \Delta_1(t)dt\label{t137}\\
&=x\Cgot \Delta_1(x)+\Cgot \Delta_1(1/x) \quad (p. p. ).\label{t138}
\end{align}
\end{prop}
\dem

D'apr\`es les estimations  \eqref{t127}, \eqref{t201}, \eqref{t202}  et le fait que 
$\Cgot$ est un op\'erateur de $L^2(0,\infty)$, les trois fonctions
\begin{align*}
x&\mapsto \Cgot \Delta_1(x)\\
x &\mapsto \big(\Cgot \Delta_1(1/x)\big)/x\\
x&\mapsto A(x)/x=A_1(x)
\end{align*}
appartiennent  à $L^2(0,\infty)$. Leurs transformées de Mellin-Plancherel respectives sont
$$
\frac{\zeta(1-s)}{\zeta(s)}M\Delta_1(1-s)=\frac{\zeta(s)\zeta(1-s)}{s} \, , \quad \frac{\zeta(s)\zeta(1-s)}{1-s} \, \text{  et  } \frac{\zeta(1-s)\zeta(s)}{s(1-s)},
$$
et \eqref{t138} équivaut à l'identité
\begin{equation*}
\frac{\zeta(s)\zeta(1-s)}{s}+\frac{\zeta(s)\zeta(1-s)}{1-s}=\frac{\zeta(1-s)\zeta(s)}{s(1-s)}.
\end{equation*}

Quant à \eqref{t137}, elle découle de l'identité
$$
\frac{\zeta(1-s)\zeta(s)}{s(1-s)}=\frac{1}{1-s}\cdot\frac{\zeta(s)\zeta(1-s)}{s},
$$
qui donne
$$
A_1(x)=\frac 1x\int_0^x \Cgot \Delta_1(t) dt \quad (p. p. )
$$
(cf. \cite{baez-duarte-all}, (12), p. 234), et l'égalité presque partout est vraie partout, par continuité.\fin

\medskip

Nous verrons au \S\ref{t126} que la proposition \ref{prop:eq-approchee-psi1} se ramène pour l'essentiel à une forme précisée de \eqref{t138}, à savoir que pour tout $x>0$,
\begin{equation}\label{t134}
xI(x,v) +I(1/x,x^2v) \vers A(x) \quad (v\vers \infty).
\end{equation}

Nous pouvons déduire de la proposition \ref{t167} une expression utile pour la fonction $A$.
\begin{prop}\label{t139}
Pour $x>0$, on a
$$
A(x)=\int_0^{\infty}\Delta(t)\frac{\sin (2\pi tx)}{\pi t^2}dt.
$$
\end{prop}
\dem
On a,
$$
\int_0^{\infty}\Delta(t)\frac{\sin (2\pi tx)}{\pi t^2}dt =\int_0^{\infty}\Delta_1(t)\frac{\sin (2\pi tx)}{\pi t}dt.
$$
Comme
$$
\frac{\sin (2\pi tx)}{\pi t}=2\int_0^x\cos(2\pi tu) \, du
$$
est la transformée de Fourier en cosinus de $u\mapsto [u\ioe x]$, le théorème de Plancherel (pour $\Cgot$ sur $L^2(0,\infty)$) nous donne
\begin{equation*}
\int_0^{\infty}\Delta(t)\frac{\sin (2\pi tx)}{\pi t^2}dt=\int_0^x\Cgot\Delta_1(t)dt.\fine
\end{equation*}

\subsection{Estimations des fonctions cosinus intégral généralisées}

Nous utiliserons des estimations de
\begin{equation*}
\ci(s,v)=\int_v^{\infty}\cos(t) \cdot t^{s-1}dt\, ; \quad \Ci(s,v)=\int_0^{v}\cos(t) \cdot t^{s-1}dt
\end{equation*}
pour $v>0$ et $s=\demi+i\tau$, $\tau \in \Real$, ce que nous supposons dans tout ce sous-paragraphe.

Notons pour commencer que $\Ci$ est absolument convergente et $\ci$ semi-convergente. De plus,

\begin{align*}
\ci(s,v)+\Ci(s,v)&=\int_0^{\infty}\cos(t) \cdot t^{s-1}dt\\
&=\Gamma(s)\cos(\pi s/2)
\end{align*}
(cf. \cite{37.0450.01}, \S 62, (4)). On a donc
\begin{align}
2(2\pi)^{-s}\big(\ci(s,v)+\Ci(s,v)\big)&=2(2\pi)^{-s}\Gamma(s)\cos(\pi s/2)\notag\\
&=\frac{\zeta(1-s)}{\zeta(s)},\label{t143}
\end{align}
d'après l'équation fonctionnelle de la fonction $\zeta$. 
\medskip

Les estimations que nous utiliserons seront démontrées grâce à la proposition classique suivante, qui découle d'une intégration par parties en écrivant $ge^{if}=(g/f')\cdot f'e^{if}$.

\begin{prop} \label{t12}
Soient $a<b$ deux nombres r\'eels,  $f,g:[a,b[ \to \Real$ deux fonctions continûment d\'erivables  telles que $g/f'$ est bien d\'efinie et monotone, et qu'il existe $c>0$ telle que pour tout $t \in[a,b]$, 
$|f'(t)/g(t)|\soe c$. On a alors
\[
\Big \vert \int_a^b g(t ) e^{if(t )} dt \Big \vert \ioe \frac{2}{c}. 
\]
\end{prop} 

La proposition suivante rassemble les estimations dont nous aurons l'usage pour les fonctions $\Ci$ et $\ci$ dans la preuve de la proposition \ref{prop:eq-approchee-psi1}. 
\begin{prop}\label{t131}
Pour $v>0$ et $\tau \in \Real$, on a
\begin{align*}
|\Ci(1/2+i\tau,v)|&\ioe \min\big(4,2\sqrt{v}/(|\tau|-v)\big) \quad (|\tau|>v)\\
|\ci(1/2+i\tau,v)|&\ioe \min\big(4,2\sqrt{v}/(v- |\tau|)\big) \quad (|\tau|<v).
\end{align*}
\end{prop}
\dem

Nous donnons la démonstration pour $\Ci$, celle pour $\ci$ étant similaire. 

On peut supposer $\tau >v>1$. On a
\begin{align}
\Ci(1/2+i\tau,v)&=\int_0^v\frac{\cos t}{\sqrt{t}}t^{i\tau}dt\notag\\
&=\int_0^v\frac{e^{i(t+\tau\log t)}}{2\sqrt{t}}dt +\int_0^v\frac{e^{i(-t+\tau\log t)}}{2\sqrt{t}}dt.\label{t17}
\end{align}

Commençons par la première intégrale. Avec les notations de la proposition \ref{t12}, on a ici
\begin{align*}
g(t) &=t^{-1/2},\\
f(t) &=t+\tau\log t,
\end{align*}
de sorte que
$$
f'(t)/g(t)= \sqrt{t}+\tau/\sqrt{t}.
$$
Cette fonction est décroissante sur $]0,\tau]$, donc sur $]0,v]$, et on a $|f'(t)/g(t)|\soe (\tau +v)/\sqrt{v}$ sur ce dernier intervalle. Par conséquent,
\begin{equation}\label{t15}
\Big \vert\int_0^v\frac{e^{i(t+\tau\log t)}}{\sqrt{t}}dt\Big \vert \ioe \frac{2\sqrt{v}}{\tau +v}.
\end{equation}
De même,
\begin{equation}\label{t132}
\Big \vert\int_0^v\frac{e^{i(-t+\tau\log t)}}{\sqrt{t}}dt\Big \vert \ioe \frac{2\sqrt{v}}{\tau -v}.
\end{equation}
En ajoutant \eqref{t15} et \eqref{t132}, on obtient 
\begin{equation}\label{t133}
|\Ci(1/2+i\tau,v)| \ioe 2\sqrt{v}/(\tau-v).
\end{equation}

Si $v$ et $\tau$ sont trop proches, on peut améliorer cette inégalité de la façon suivante. Si $\tau \soe v+\sqrt{v}$, \eqref{t133} nous donne
$$
|\Ci(1/2+i\tau,v)| \ioe 2.
$$
En revanche, si $\tau-\sqrt{v}<v<\tau$, on a
\begin{align*}
|\Ci(1/2+i\tau,v)| &\ioe |\Ci(1/2+i\tau,\tau-\sqrt{v})| +\int_{\tau-\sqrt{v}}^v\frac{dt}{\sqrt{t}}\\
&\ioe 2+\int_{v-\sqrt{v}}^v\frac{dt}{\sqrt{t}}\\
&\ioe 4.\fine
\end{align*}

\subsection{Démonstration de la proposition \ref{prop:eq-approchee-psi1}}\label{t126}

Pour $x>0$, $v>0$, nous posons
$$
\psi_1(x,v)=-\frac{1}{\pi}\sum_{n\le  v} \frac{\tau(n)}{n} \sin(2\pi nx)
$$
(sommes partielles de la série $\psi_1(x)$). La proposition suivante ramène l'étude de $\psi_1(x,v)$ à celle de l'intégrale $I(x,v)$ définie par \eqref{t136}

\begin{prop}\label{t140}
Pour $x>0$, $v>0$, on a
\begin{equation}\label{t142}
-\psi_1(x,v)=-xI(x,v)+A(x)-\demi(\log x +\log 2\pi -\gamma)+\eps(x,v),
\end{equation}
où
$$
\eps(x,v)=-\int_v^{\infty}\frac{\sin (2\pi tx)}{\pi t}(\log t +2\gamma) dt + \frac{\sin (2\pi vx)}{\pi v}\Delta(v)-\int_v^{\infty}\Delta(t)\frac{\sin (2\pi tx)}{\pi t^2}dt.
$$
\end{prop}
\dem

Nous commençons par effectuer une intégration par parties sur l'expression de $-\psi_1(x,v)$ comme intégrale de Stieltjes :
\begin{align}
-\psi_1(x,v)&=\frac{1}{\pi}\sum_{n\le  v} \frac{\tau(n)}{n} \sin(2\pi n x)\notag\\
&=\int_0^v\frac{\sin (2\pi tx)}{\pi t}d\big (\sum_{n\ioe t}\tau(n)\big )\notag\\
&=\int_0^v\frac{\sin (2\pi tx)}{\pi t}(\log t +2\gamma) dt+\int_0^v\frac{\sin (2\pi tx)}{\pi t}d\Delta(t).\label{t6}
\end{align}

Pour la première intégrale de \eqref{t6}, on a
$$
\int_0^v\frac{\sin (2\pi tx)}{\pi t}(\log t +2\gamma) dt=\int_0^{\infty}-\int_v^{\infty}.
$$
D'une part, nous posons
$$
\eps_1(x,v)=-\int_v^{\infty}\frac{\sin (2\pi tx)}{\pi t}(\log t +2\gamma) dt
$$
et d'autre part
\begin{align}
\int_0^{\infty}\frac{\sin (2\pi tx)}{\pi t}(\log t +2\gamma) dt&=\frac{1}{\pi}\int_0^{\infty}(\log (t/2\pi x) +2\gamma) \frac{\sin t}{t} dt\notag\\
&=\frac{1}{\pi}\int_0^{\infty}\log t \frac{\sin t}{t} dt -\frac{(\log x +\log 2\pi -2\gamma)}{\pi}\int_0^{\infty}\frac{\sin t}{t} dt\notag\\
&=-\demi(\log x +\log 2\pi -\gamma),\label{t141}
\end{align}
où l'on a utilisé les valeurs de l'intégrale d'Arndt 
$$
\int_0^{\infty}\log t \frac{\sin t}{t} dt=-\frac{\pi\gamma}{2}
$$
(cf. \cite{37.0450.01} \S 69, (9)) et de celle de Dirichlet
$$
\int_0^{\infty} \frac{\sin t}{t} dt=\frac{\pi}{2}.
$$

\smallskip

Passons à la seconde intégrale de \eqref{t6} :
\begin{align}
\int_0^v\frac{\sin (2\pi tx)}{\pi t}d\Delta(t)&=\frac{\sin (2\pi tx)}{\pi t}\Delta(t) \Big\vert_0^v-\int_0^v\Delta(t)\Big (-\frac{\sin (2\pi tx)}{\pi t^2}+2x\frac{\cos (2\pi tx)}{t}\Big)dt\notag\\
&=\int_0^{\infty}\Delta(t)\frac{\sin (2\pi tx)}{\pi t^2}dt -xI(x,v)+\eps_2(x,v),\label{t9}
\end{align}
où l'on a posé
\begin{equation*}
\eps_2(x,v)=\frac{\sin (2\pi vx)}{\pi v}\Delta(v)-\int_v^{\infty}\Delta(t)\frac{\sin (2\pi tx)}{\pi t^2}dt.
\end{equation*}

En ajoutant \eqref{t141} et \eqref{t9}, et en tenant compte de la proposition \ref{t139}, on obtient bien \eqref{t142}.\fin

\medskip

La proposition suivante précise l'assertion \eqref{t134}.
\begin{prop}
Pour $x>0$, $v>0$, on a
\begin{equation}\label{t144}
xI(x,v)+I(1/x,x^2v)=A(x)+\eta(x,v),
\end{equation}
où
\begin{multline}\label{t147}
\eta(x,v)=-\int_{|\tau|\ioe V}\frac{\zeta^2(1-s)}{1-s}(2\pi )^{s-1}\ci(1-s,V)(x^s+x^{1-s})\,\frac{d\tau}{\pi}+\\\int_{|\tau|> V}\frac{\zeta^2(1-s)}{1-s}\Ci(1-s,2\pi xv)(2\pi )^{s-1}(x^s+x^{1-s})\,\frac{d\tau}{\pi}-\int_{|\tau|>V}\frac{\zeta(s)\zeta(1-s)}{s(1-s)}x^s\frac{d\tau}{2\pi},
\end{multline}
avec $V=2\pi xv$.
\end{prop}
\dem

En utilisant l'identité \eqref{t5}, on a
\begin{align*}
xI(x,v)+I(1/x,x^2v) &= \int_{\Re s=\demi}\frac{\zeta^2(1-s)}{1-s}\Ci(1-s,2\pi xv)(2\pi )^{s-1}(x^s+x^{1-s})\,\frac{d\tau}{\pi}\\
&=\int_{|\tau|\ioe V}+\int_{|\tau|> V}\, ,
\end{align*}
où $V$ est un paramètre positif arbitraire ; cela étant, la proposition \ref{t131} nous incite à poser sans attendre $V=2\pi xv$.

L'intégrale $\int_{|\tau|> V}$ est la deuxième intervenant dans $\eta(x,v)$. Quant à l'intégrale $\int_{|\tau|\ioe V}$, elle vaut d'après~\eqref{t143}
$$
\int_{|\tau|\ioe V}\frac{\zeta^2(1-s)}{1-s}\Big (\frac{\zeta(s)}{\zeta(1-s)}-2(2\pi )^{s-1}\ci(1-s,V)\Big )(x^s+x^{1-s})\,\frac{d\tau}{2\pi}=J_1-J_2,
$$
où
\begin{align*}
J_1 &=\int_{|\tau|\ioe V}\frac{\zeta(s)\zeta(1-s)}{1-s}(x^s+x^{1-s})\,\frac{d\tau}{2\pi}\, ;\\
J_2&=\int_{|\tau|\ioe V}\frac{\zeta^2(1-s)}{1-s}(2\pi )^{s-1}\ci(1-s,V)(x^s+x^{1-s})\,\frac{d\tau}{\pi}.
\end{align*}

On a 
\begin{align*}
J_1 &=\int_{|\tau|\ioe V}\frac{\zeta(s)\zeta(1-s)}{1-s}x^s\frac{d\tau}{2\pi}+  \int_{|\tau|\ioe V}\frac{\zeta(s)\zeta(1-s)}{s}x^{s}\frac{d\tau}{2\pi}\\ 
&=\int_{|\tau|\ioe V}\frac{\zeta(s)\zeta(1-s)}{s(1-s)}x^s\frac{d\tau}{2\pi}\\
&=A(x)-\int_{|\tau|>V}\frac{\zeta(s)\zeta(1-s)}{s(1-s)}x^s\frac{d\tau}{2\pi},
\end{align*}
d'après \eqref{t145}. Cela donne bien \eqref{t144}.\fin

\medskip

En appliquant \eqref{t142} aux couples $(x,v)$ et $(1/x,x^2v)$, on obtient
\begin{align*}
-\psi_1(x,v)&-x\psi_1(1/x,x^2v)\\&=-xI(x,v)+A(x)-\demi(\log x +\log 2\pi -\gamma )+\eps(x,v)\\
 & \quad +x\Big (-x^{-1}I(1/x,x^2v)+A(1/x)-\demi (\log (1/x) +\log 2\pi -\gamma )+\eps(1/x,x^2v)\Big)\\
&=-xI(x,v)-I(1/x,x^2v)+2A(x)+\frac{x-1}{2}\log x-\frac{x+1}{2}\big(\log(2\pi) -\gamma\big)+\\
&\quad +\eps(x,v)+x\eps(1/x,x^2v)\\
&=-\demi\log x -F(x) +\eps(x,v)+x\eps(1/x,x^2v)-\eta(x,v),
\end{align*}
où l'on a utilisé \eqref{t144} et la définition \eqref{eq:def-F} de la fonction $F$.

Pour conclure la démonstration de la proposition \ref{prop:eq-approchee-psi1}, il nous reste à estimer les fonctions $\eps$ et $\eta$. C'est l'objet des deux propositions suivantes.
\begin{prop}\label{t146}
Pour $0<x<1$ et $x^2v\soe 2$, on a
$$
\eps(x,v)+x\eps(1/x,x^2v) \ll (x^2v)^{-1/2}.
$$
\end{prop}
\dem

Nous commençons par estimer $\eps(x,v)$ sous la seule hypothèse $v\soe 2$. Rappelons que
$$
\eps(x,v)=-\int_v^{\infty}\frac{\sin (2\pi tx)}{\pi t}(\log t +2\gamma) dt + \frac{\sin (2\pi vx)}{\pi v}\Delta(v)-\int_v^{\infty}\Delta(t)\frac{\sin (2\pi tx)}{\pi t^2}dt.
$$
Par le second théorème de la moyenne, la première intégrale est $\ll (\log v)/xv$. En utilisant l'estimation de Dirichlet $\Delta(v) \ll v^{1/2}$, on voit que les autres termes sont $\ll v^{-1/2}$.

On en déduit, si $0<x<1$ et $x^2v\soe 2$ :
\begin{align*}
\eps(x,v)+x\eps(1/x,x^2v) &\ll (\log v)/xv +v^{-1/2} +x\log(x^2v)/xv +x(x^2v)^{-1/2}\\
&\ll (x^2v)^{-1/2}.\fine
\end{align*}

\begin{prop}
Pour $0<x<1$ et $x^2v\soe 2$, on a
$$
\eta(x,v) \ll (x^2v)^{-1/2}\log^2(x^2v).
$$
\end{prop}
\dem

Pour $0<x<1$ et $x^2v\soe 2$, on a $V=2\pi xv>12$. En utilisant la proposition \ref{t131}, on voit sur la définition \eqref{t147} de $\eta(x,v)$ que
\begin{equation}\label{t148}
\eta(x,v) \ll x^{1/2}\int_{0}^{\infty}\frac{|\zeta(\demi+i\tau)|^2}{1+\tau}\min(1,V^{1/2}\big\vert V-\tau\big\vert^{-1})d\tau.
\end{equation}

Une estimation du type  
$$
\zeta\left(\demi+i\tau\right)\ll (1+|\tau|)^{\delta}
$$
(cf. par exemple le corollaire 3.7, p. 234 de \cite{MR1366197}) conduit à
$$
\int_{0}^{\infty}\frac{|\zeta(\demi+i\tau)|^2}{1+\tau}\min(1,V^{1/2}\big\vert V-\tau\big\vert^{-1})d\tau \ll V^{2\delta-1/2}\log V,
$$
ce qui fournit un terme d'erreur admissible dans \eqref{eq:fonctionnelle-approche} dès que $\delta <1/4$. Cependant, on aboutit à un meilleur résultat en utilisant une estimation des fonctions $I_2(T)$ et $E(T)$ définies par 
$$
I_2(T)=\int_0^T|\zeta(1/2+i\tau)|^2d\tau=T\log T-(\log 2\pi +1-2\gamma)T+E(T) \quad (T>0).
$$

La fonction $E(T)$ est l'objet d'une abondante littérature (cf. par exemple \cite{ivic}, chapters 4, 15) ; nous nous contenterons de l'estimation classique d'Ingham : $E(T)\ll T^{1/2}\log T$ pour $T\soe 2$ (cf. \cite{53.0313.01}, Theorem A', p. 294).

La contribution de l'intervalle $|\tau -V|\ioe \sqrt{V}$ à l'intégrale de \eqref{t148} est 
\begin{align*}
 &\ll V^{-1}\int_{V-\sqrt{V}}^{V+\sqrt{V}}|\zeta(1/2+i\tau)|^2d\tau\\
&\ll V^{-1/2}\log V,
\end{align*}
d'après l'estimation d'Ingham.

La contribution des $\tau \soe 2V$ est 
\begin{align*}
 &\ll V^{1/2}\int_{2V}^{\infty}\frac{|\zeta(1/2+i\tau)|^2}{\tau^2}   d\tau\\
&\ll V^{-1/2}\log V,
\end{align*}
en utilisant simplement l'estimation $I_2(T)\ll T\log T$ pour $T\soe 2$ et une intégration par parties. 

De même la contribution des $\tau \ioe V/2$ est 
\begin{align*}
 &\ll V^{-1/2}\int_0^{V/2}\frac{|\zeta(1/2+i\tau)|^2}{1+\tau}  d\tau\\
&\ll V^{-1/2}\log^2 V.
\end{align*}

Maintenant, la contribution de l'intervalle $V+\sqrt{V}<\tau<2V$ est 
\begin{align*}
 &\ll V^{-1/2}\int_{V+\sqrt{V}}^{2V}\frac{|\zeta(1/2+i\tau)|^2}{\tau -V}   d\tau\\
&\ioe V^{-1/2}\sum_{1\ioe k\ioe \sqrt{V}}\frac{1}{k\sqrt{V}}\int_{V+k\sqrt{V}}^{V+(k+1)\sqrt{V}}|\zeta(1/2+i\tau)|^2 d\tau\\
&\ll V^{-1/2}(\log V)\sum_{1\ioe k\ioe \sqrt{V}}\frac 1k\expli{d'après l'estimation d'Ingham}\\
&\ll V^{-1/2}\log^2 V,
\end{align*}
et on a encore la même estimation pour la contribution de l'intervalle $V/2<\tau<V-\sqrt{V}$. 

En conclusion, on a pour $0<x\ioe 1$ et $x^2v\soe 2$,
\begin{align*}
\eta(x,v) &\ll x^{1/2} \cdot \frac{\log^2(xv)}{\sqrt{xv}}\\
&\ll \frac{\log^2(x^2v)}{\sqrt{x^2v}}.\fine
\end{align*}

\section{Conclusion de la preuve du th\'eor\`eme \ref{t150}}\label{t149}

D'après la proposition \ref{prop:lien-fhi1-W}, la série $\fhi_1(x)$ converge si seulement la série 
$\W(x)$ converge, et dans ce cas, 
$$
\fhi_1(x)=-\frac12 \W(x)+G(x)+\delta(x). 
$$
Par ailleurs la somme partielle $\psi_1(x,v)$ de la s\'erie $\psi_1(x)$ satisfait aux points i), ii) et iii) \'enonc\'es dans le paragraphe \ref{par:eq-approchee}. 
En effet, le point i) avec les param\`etres $a=2$, $b=1/2-\eps$ (quel que soit $\eps>0$) résulte de la proposition \ref{prop:eq-approchee-psi1}. Le point ii) est trivial. En ce qui concerne le point iii), on dispose de la majoration  de Walfisz  (cf. \cite{walfisz}, $(25_{VI})$, p. 566)
\begin{equation*}
\sum_{n\le v} \frac{\tau(n)}{n} \sin(2\pi n x) =\sum_{m\le v} \frac{1}{m}\sum_{k\ioe v}\frac{1}{k} \sin(2\pi km x)\ll \log (v) 
\end{equation*}
uniformément pour $x \in\Real$ et $v\ge 2$.

On obtient donc, d'apr\`es la proposition \ref{prop:meta-th}, que les séries $\psi_1(x)$ et $\W(x)$ convergent simultanément et qu'en cas de convergence  on a  
\[
\psi_1(x)=-\frac12 \W(x)+G(x)+\delta(x). 
\]
Cela achève la preuve du théorème \ref{t150}. 

\bigskip

\begin{center}
  {\sc Remerciements}
\end{center}
\begin{quote}
{\footnotesize Outre leurs laboratoires respectifs, les auteurs remercient les institutions ayant favorisé leur travail sur cet article : le laboratoire franco-russe Poncelet (CNRS, Université Indépendante de Moscou), et l'Institut Mittag-Leffler (Djursholm).
}
\end{quote}

\medskip

\begin{multicols}{2}
\footnotesize

\noindent BALAZARD, Michel\\
Institut de Math\'ematiques de Luminy\\
CNRS, Universit\'e de la M\'editerran\'ee\\
Campus de Luminy, Case 907\\
13288 Marseille Cedex 9\\
FRANCE\\
Adresse \'electronique : \texttt{balazard@iml.univ-mrs.fr}

\smallskip

\noindent MARTIN, Bruno\\
Laboratoire de Math\'ematiques Pures et Appliqu\'ees\\
CNRS, Universit\'e du Littoral C\^ote d'Opale\\
50 rue F. Buisson, BP 599\\
62228 Calais Cedex\\
FRANCE\\
Adresse \'electronique : \texttt{martin@lmpa.univ-littoral.fr}
\end{multicols}

\end{document}